%%%%%%%%%%% pi-two.tex %%%% revised version 4 %%% 17/10/04

%%%%%%%%%%%%%%%%%%%%%%%%%%%%%%%%%%%%%%%%%%%%%%%%%%%%%%%%%%%%%%%
%%%%%%%%%%%             gtmacros.tex            %%%%%%%%%%%%%%%
%%%%%%%%%%%             version 1.7             %%%%%%%%%%%%%%% 
%
%                       Colin Rourke   
%
%
%    These macros are recommended for use by authors submitting articles   
%    to Geometry and Topology or to Algebraic and Geometric Topology.  
%    They are intended to be used with plain TeX. Each macro is described 
%    briefly to make it clear how to use it (or to modify it to achieve
%    different results).  If you modify this file then please change its
%    name.  If you modify this file and use the modified file to 
%    format an article for submission to Geometry and Topology or
%    Algebraic and Geometric Topology, then please paste the modified
%    file into your main TeX file.  Do not submit it as a separate file.
%      
%    Instructions on using these macros are also given in  gtmacins.tex  
%    or  gtmacins.ps  or .pdf  available on the gt www pages or by 
%    anonymous ftp from the gt/info/macros directory.
%
%
\magnification=\magstephalf      % Sets default point size to 11pt.
%
%  Basic layout parameters :
%
\vsize=7.5truein                 % Sets text height to 7.5 inches.
\hsize=5.2truein                 % Sets text width to 5.2 inches.
\newskip\stdskip                 % standard vertical space
\stdskip=6pt plus3pt minus3pt    % (slightly more stretchy
\medskipamount=\stdskip          % than the usual \medskip)
\parindent=0pt                   % Paragraphs are non-indented with
\parskip=\stdskip                % a little space between paragraphs. 
\abovedisplayskip=\stdskip       %  Reduces the space
\belowdisplayskip=\stdskip       %  around displays.
\mathsurround=0.75pt             % Gives a little extra space around maths.
\overfullrule=0pt                %  Prevents black boxes
%
%   The following macro is for principal paragraph breaks ie
%   a paragraph break with a slightly larger space :
%
\def\ppar{\par\goodbreak\vskip 8pt plus 4pt minus 4pt}     
%
%  The standard horizontal space for theorems, labels etc :
%
\def\stdspace{\hskip 0.75em plus 0.15em\ignorespaces}
\let\qua\stdspace % useful abbreviation (3/4 of a quad)
%
%%%%%%%%%%%%%%            FONT MACROS            %%%%%%%%%%%%
%
%           The following font macros define the AMS symbol 
%           and Euler-Fraktal fonts for use in text and
%           mathematics with appropriate size changes.
%           They also define two new control sequences  
%           \small  and  \large  (similar to those built
%           into LaTeX) which change the size of all fonts 
%           both in text and maths.  \small  is 10% smaller 
%           than normal and  \large  30% bigger.  The strange
%           size of the \large text fonts (10pt scaled 1315)
%           is because these macros are intended to be used
%           at \magstephalf.  The result is 10pt scaled 1440
%           (\magstep2) which is a standard font size.  If
%           you are borrowing these macros to use them at
%           another basic  \magnification, then you will
%           probably need to change 1315 to 1200 in the eleven
%           places marked ** below.  \large  will then be
%           20% bigger than normal.  Note that at \magstephalf
%           all the fonts come out roughly one point larger
%           than their size as defined in these macros.
%
%           The size-changing macros are based on Knuth's
%           \ninepoint and \eightpoint macros.
%
%
%    The macros are laid out in a way which makes it clear how to
%    add futher fonts (or delete unavailable ones) and how to add
%    further size changes.
%
%    First comes a definition of  \hexnumber  which is needed to
%    refer to font families whose family number is not known :
%
\def\hexnumber#1{\ifcase#1 0\or 1\or 2\or 3\or 4\or 5\or 6\or 7\or 8\or
 9\or A\or B\or C\or D\or E\or F\fi}
%
%     Next we define the AMS symbol-a fonts at 13,10,9,7,6,5 points
%
\font\thirtnmsa=msam10 scaled 1315    %%% **  see note above 
\font\tenmsa=msam10          \font\ninemsa=msam9
\font\sevenmsa=msam7         \font\sixmsa=msam6
\font\fivemsa=msam5
%%%%%%  (add further sizes here if you need them)
%
%    and the standard size family for these fonts
%
\newfam\msafam                  \textfont\msafam=\tenmsa
\scriptfont\msafam=\sevenmsa    \scriptscriptfont\msafam=\fivemsa
\edef\hexa{\hexnumber\msafam}        %  The msa family is  \fam\hexa
\def\msa{\fam\msafam\tenmsa}         %  \msa  switches to this family
%
%    Repeat these steps for the AMS symbol-b fonts
%
\font\thirtnmsb=msbm10 scaled 1315   %%%  ** see note above
\font\tenmsb=msbm10      \font\ninemsb=msbm9
\font\sevenmsb=msbm7     \font\sixmsb=msbm6
\font\fivemsb=msbm5
%%%%%%  (add further sizes here if you need them)
%
\newfam\msbfam                   \textfont\msbfam=\tenmsb       
\scriptfont\msbfam=\sevenmsb     \scriptscriptfont\msbfam=\fivemsb
\edef\hexb{\hexnumber\msbfam}    %  The msb family is \fam\hexb  
\def\msb{\fam\msbfam\tenmsb}     %  \msb switches to this family
%
%        Repeat for the Euler-Fraktal fonts 
%
\font\thirtneufm=eufm10 scaled 1315   %%% **  see note above 
\font\teneufm=eufm10                 \font\nineeufm=eufm9
\font\seveneufm=eufm7                \font\sixeufm=eufm6
\font\fiveeufm=eufm5
%%%%%%  (add further sizes here if you need them)
%
\newfam\eufmfam                    \textfont\eufmfam=\teneufm
\scriptfont\eufmfam=\seveneufm     \scriptscriptfont\eufmfam=\fiveeufm
\edef\hexf{\hexnumber\eufmfam}      % The Euler-Fraktal family is
\def\frak{\fam\eufmfam\teneufm}     % \fam\hexf and \frak switches to this
%
%%%  Add further fonts families here (using the same format) if you need
%    them.  The def of hexnumber is optional (it is only used for
%    \mathchardef 's).
%
%      Now we need to define the standard fonts (which are
%      already defined at 10,7 and 5 point) at 13,9 and 6 point:
%
%      Roman fonts:
\font\thirtnrm=cmr10 scaled 1315    %%%  ** see note above
\font\ninerm=cmr9                   \font\sixrm=cmr6   
%%%%%%  (add further sizes here if you need them)
%
%      Math italic fonts
\font\thirtni=cmmi10 scaled 1315    %%%  ** see note above 
\font\ninei=cmmi9                   \font\sixi=cmmi6  
%%%%%%  (add further sizes here if you need them)
%
%     Symbol fonts
\font\thirtnsy=cmsy10 scaled 1315   %%%  ** see note above
\font\ninesy=cmsy9                  \font\sixsy=cmsy6  
%%%%%%  (add further sizes here if you need them)
%
%     Bold face
\font\thirtnbf=cmbx10 scaled 1315   %%%  ** see note above 
\font\ninebf=cmbx9                  \font\sixbf=cmbx6  
%%%%%%  (add further sizes here if you need them)
%
%     The maths extension font (only defined at text size)
%
\font\thirtnex=cmex10 scaled 1315   %%%  ** see note above
\font\nineex=cmex9                  
%%%%%%  (add further sizes here if you need them)
%
%     Finally three fonts (text italic, slanted and typewriter type)
%     which are also only defined at text size
%
\font\thirtnit=cmti10 scaled 1315  %%%  ** see note above 
\font\nineit=cmti9                  
%%%%%%  (add further sizes here if you need them)
%
\font\thirtnsl=cmsl10 scaled 1315  %%%  ** see note above 
\font\ninesl=cmsl9                  
%%%%%%  (add further sizes here if you need them)
%
\font\thirtntt=cmtt10 scaled 1315  %%%  ** see note above 
\font\ninett=cmtt9                  
%%%%%%  (add further sizes here if you need them)
%
%
%     Now come the two main macros.  What  \small  does is to
%     change all the families of fonts from normal size which is
%     10,7,5  (ie 10pt text, 7pt subscript, 5pt subsubscript)
%     to 9,6,5.  \large  similarly changes to  13,9,7.  To make
%     other size changing macros, choose your three sizes, add
%     font size definitions if necessary and make the obvious changes
%     to one of these macros.  Change  \normalbaselineskip  and
%     \strutbox  dimensions to appropriate sizes as well.  To
%     add further fonts, insert them in each macro, using the
%     AMS fonts as a model.
%      
%
\def\small{%
%
%   redefine the sizes of the roman fonts :
%
\textfont0=\ninerm \scriptfont0=\sixrm \scriptscriptfont0=\fiverm
\def\rm{\fam0\ninerm}%       % ( \rm  sets \ninerm  in text mode
%                            %  and \fam0 in math mode)
%
%   and the math italic fonts :
%
\textfont1=\ninei \scriptfont1=\sixi \scriptscriptfont1=\fivei
%
%   and the symbol fonts :
%
\textfont2=\ninesy \scriptfont2=\sixsy \scriptscriptfont2=\fivesy
%
%   There is only one math extension font :
%
\textfont3=\nineex \scriptfont3=\nineex \scriptscriptfont3=\nineex
%
%   Next the bold font (named rather than numbered) :
%
\textfont\bffam=\ninebf \scriptfont\bffam=\sixbf
\scriptscriptfont\bffam=\fivebf \def\bf{\fam\bffam\ninebf}%
%
%   and the three text-only fonts : 
%
\textfont\itfam=\nineit \def\it{\fam\itfam\nineit}%
\textfont\slfam=\ninesl \def\sl{\fam\slfam\ninesl}%
\textfont\ttfam=\ninett \def\tt{\fam\ttfam\ninett}%
%
%   Now the three new families of AMS fonts :
%
%   AMS symbol-a
%
\textfont\msafam=\ninemsa \scriptfont\msafam=\sixmsa
\scriptscriptfont\msafam=\fivemsa \def\msa{\fam\msafam\ninemsa}%         
%
%   AMS symbol-b
%
\textfont\msbfam=\ninemsb \scriptfont\msbfam=\sixmsb
\scriptscriptfont\msbfam=\fivemsb \def\msb{\fam\msbfam\ninemsb}%         
%
%   Euler-Fraktal font
%
\textfont\eufmfam=\nineeufm  \scriptfont\eufmfam=\sixeufm
\scriptscriptfont\eufmfam=\fiveeufm \def\frak{\fam\eufmfam\nineeufm}%
%
%%%  Add further fonts families here if you need them.
%
%    Reset \normalbaselineskip and \strubox :
%
\normalbaselineskip=11pt%
\setbox\strutbox=\hbox{\vrule height8pt depth3pt width0pt}%
%
%    Set \normalbaselines and \rm (roman) as defaults :
%
\normalbaselines\rm
%
%    Reset some of the basic vertical skips:
%
\stdskip=4pt plus2pt minus2pt    
\medskipamount=\stdskip          
\parskip=\stdskip                
\abovedisplayskip=\stdskip       
\belowdisplayskip=\stdskip       
\def\ppar{\par\goodbreak\vskip 6pt plus 3pt minus 3pt}%     
%
%   And finally reset the size of section heads (see below):
%
\def\section##1{\global\advance\sectionnumber by 1
\vskip-\lastskip\penalty-800\vskip 20pt plus10pt minus5pt 
\egroup{\bf\number\sectionnumber\quad##1}\bgroup\small         
\vskip 6pt plus3pt minus3pt
\nobreak\resultnumber=1}%      % Reset resultnumber at start of section
}    %%%   End of  \small  macro      
%
%   Two useful abbreviations to keep track of \small material:
\def\beginsmall{\bgroup\small}
\let\endsmall\egroup
%
%
%    The \large  macro is similar (comments abbreviated):
%
%
\def\large{%
\textfont0=\thirtnrm \scriptfont0=\ninerm \scriptscriptfont0=\sevenrm
\def\rm{\fam0\thirtnrm}%
\textfont1=\thirtni \scriptfont1=\ninei \scriptscriptfont1=\seveni
\textfont2=\thirtnsy \scriptfont2=\ninesy \scriptscriptfont2=\sevensy
\textfont3=\thirtnex \scriptfont3=\thirtnex \scriptscriptfont3=\thirtnex
\textfont\bffam=\thirtnbf \scriptfont\bffam=\ninebf
\scriptscriptfont\bffam=\sevenbf \def\bf{\fam\bffam\thirtnbf}%
\textfont\itfam=\thirtnit \def\it{\fam\itfam\thirtnit}%
\textfont\slfam=\thirtnsl \def\sl{\fam\slfam\thirtnsl}%
\textfont\ttfam=\thirtntt \def\tt{\fam\ttfam\thirtntt}%
%   AMS symbol-a  :
\textfont\msafam=\thirtnmsa \scriptfont\msafam=\ninemsa
\scriptscriptfont\msafam=\sevenmsa \def\msa{\fam\msafam\thirtnmsa}%         
%   AMS symbol-b  :
\textfont\msbfam=\thirtnmsb \scriptfont\msbfam=\ninemsb
\scriptscriptfont\msbfam=\sevenmsb \def\msb{\fam\msbfam\thirtnmsb}%         
%   Euler-Fraktal font :
\textfont\eufmfam=\thirtneufm  \scriptfont\eufmfam=\nineeufm
\scriptscriptfont\eufmfam=\seveneufm \def\frak{\fam\eufmfam\teneufm}%
%%%% Add further fonts families here if you need them.
%   Reset \normalbaselineskip and \strubox and initialise :
\normalbaselineskip=16pt%
\setbox\strutbox=\hbox{\vrule height11.5pt depth4.5pt width0pt}%
\normalbaselines\rm}%
%     
   %  for compatibility with latex
%
%   The next two lines define commonly used switches for
%   blackboard bold (\Bbb) and gothic type (\goth).  The   
%   \Bbb  switch is set to work in the same way as in amstex
%   and switches only the next character to blackboard bold.
\def\Bbb#1{{\msb#1}}

%
%   To use the new AMS fonts you can either use the control
%   sequences \msa, \msb (alias \Bbb) and \frak (alias \goth) eg :

   % see the msam font table
%
%   or, more generally, make \mathchardef's (cf Knuth p155) eg :
\mathchardef\plussquare="0\hexa01
\mathchardef\nge="3\hexb0B
\mathchardef\maltesecross="0\hexa7A
\mathchardef\del="0\hexf01
%
%   or you can use the amstex names for all the new symbols by
%   inserting the line  \input amsnames  in your file directly
%   after \input gtmacros. 
%   This presupposes that you have collected a copy of the file
%   amsnames.tex  from the  gt/info/macros  ftp directory.
%
%
%   Finally we need a small capital font (for author(s)) :
%
\font\sc=cmcsc10
%
%%%%%%%%%%%%%%%%%       END OF FONT MACROS     %%%%%%%%%%%%%
%
%
%                 Knuth's \square macro :
%
\def\sqr#1#2{{\vcenter{\vbox{\hrule  height.#2truept
	\hbox{\vrule width.#2truept height#1truept 
	\kern#1truept \vrule width.#2truept}
	\hrule height.#2truept}}}}
\def\sq{\sqr55}    %   A small square for end-of-proofs. 
%                  %   (Define other size squares by varing the
%                  %   the two numbers.)
%
%
%      Style macros for section heads, theorem statements etc :
%   
%
\newcount\sectionnumber            %%%  Allocate registers to take
\newcount\resultnumber             %%%  section and result numbers.
\sectionnumber=0\resultnumber=1    %%%  Set these registers to 0 and 1
%
%   The \section macro produces a \large bold faced section heading
%   numbered to the left.  Pagebreaks are encouraged before the
%   start of the section and discouraged directly after the heading.
%   Typical use  \section{First steps}  with typical result :
%
%    1  First Steps     (set bold and \large)
%
\def\section#1{\global\advance\sectionnumber by 1
\xdef\nextkey{\number\sectionnumber}%      (used by the \key macro)
\vskip-\lastskip\penalty-800\vskip 20pt plus10pt minus5pt 
{\large\bf\number\sectionnumber\quad#1}         
\vskip 8pt plus4pt minus4pt
\nobreak\resultnumber=1}      % Reset resultnumber at start of section
%
%
%
%   Next a macro to set subheadings (like the  \section  macro
%   but without the number, with less space and set in standard size).
%
%   Typical use :  \sh{Example formats}
%
         
%
%   The \proc ... \endproc macros ("proclaim") are for setting theorems, 
%   lemmas, conjectures etc with automatic numbering.  Typical use :    
%  
%    \proc{Theorem}Every lemon is yellow.\endproc
%
%   Typical result :
%     
%    Theorem 3.4  Every lemon is yellow.   

%   (with Theorem 3.4 set bold and a \stdspace of space before the 
%   statement set in slanted type).
%
\def\proc#1{\xdef\nextkey{\number\sectionnumber.\number\resultnumber}%
\vskip-\lastskip\ppar\bf%
\noindent#1\ \number\sectionnumber.\number\resultnumber
\stdspace\sl\global\advance\resultnumber by 1\ignorespaces}
\def\endproc{\rm\ppar} 
%
%  The \prf ... \endprf macros are for setting proofs.  The code
%  for \prf includes the code for \endproc, so there is no need to
%  type \endproc if the theorem is followed immediatedly by a proof.
%
\def\prf{\vskip-\lastskip\ppar\noindent{\bf Proof}%
\stdspace\rm}                            %  For start of proofs  
\def\qed{\hfill$\sq$\par\goodbreak\rm}   %  For end (or absence) of proofs
\def\endprf{\unskip\stdspace\hbox{}%     %  For end of proof (with
\hfill$\sq$\par\medskip}                 %  extra vertical space)  
        %  For start of proof with alternative name
              %  \endproof is an alias for \endprf
%
%   Typical uses :    
%  
%    \proc{Theorem}Every lemon is yellow. \qed\endproc
%
%    \proc{Theorem}Every lemon is yellow.
%    \prf Use your eyes. \endprf
%
%    \proc{Theorem}Every lemon is yellow.
%    \proof{Proof of theorem} Use your eyes. \endprf
%
%   The next macro is a variant of the \proc macro.  It has
%   exactly the same result except that it omits the number.
%
%   Typical use :  
%    
%    \proclaim{Conjecture}Some oranges are yellow.\endproc
%
\def\proclaim#1{\vskip-\lastskip\ppar\bf%
\noindent#1\stdspace\sl\ignorespaces} 

%
%   The next macro is a further variant for remarks, definitions etc.   
%   It omits the number and does not switch on slanted type.  
%  
%   Typical use :
%
%    \rk{Remark}Some lemons are thick-skinned.\endrk
%
\def\rk#1{\vskip-\lastskip\ppar{\bf #1}\stdspace\ignorespaces}                
\def\endrk{\par\medskip}
%
%   The next macro is for numbering equations etc, \label  produces the 
%   correct number  x.y  and advances the resultnumber register
%
%   Typical use :
%
%     $$fx=7\eqno{\bf\label}$$
%
%   result :
%
%                           fx = 7                           3.5
%
\def\label{\xdef\nextkey{\number\sectionnumber.\number\resultnumber}%
\number\sectionnumber.\number\resultnumber
\global\advance\resultnumber by 1}
%
%
%
%   The next macros are to automate external references.  To use them 
%   type \reflist ..... \endreflist near the beginning of your paper, 
%   where  .... is the list of references in alphabetical order 
%   and in  the form  \key{KEY}  reference    where "KEY" is a 
%   string of characters which reminds you of the reference.   
%   Separate  references with a blank line or a \par.   Eg 
%
%     \reflist
%
%     ..... more references ....
%
%     \key{Kn-84} {\bf D Knuth}, {\it The TeXbook}, Addison--Wesley (1984)
%
%     ..... more references ....
%
%     \endreflist
%
%   Then type  \references  where you wish the references to be printed
%   (normally near the end of the paper).  To refer to Knuth type
%   for example    see Knuth [\ref{Kn-84}, page 320]   and the correct
%   numerical reference will be printed.  Edit the \references macro
%   to change the formatting (if desired).
%   There is an alternative \refkey for \key, provided your KEY contains
%   only letters.  The syntax is:
%
%     \reflist
%
%     ..... more references ....
%
%     \refkey\Knuth  {\bf D Knuth}, {\it The TeXbook}, Addison--Wesley (1984)
%
%     ..... more references ....
%
%     \endreflist
%
%   \key{Knuth}  has exactly the same maening as \refkey\Knuth and you
%   can mix the two syntaxes if you want.  But \refkey\Kn-84
%   would not work.  It would set Kn as the KEY and -84 would get printed!
%
\newcount\refnumber              %  Register for reference numbers
\refnumber=1                     %  set initially to 1.
\long\def\reflist#1\endreflist{%
\long\def\thereflist{#1}{\def\refkey##1##2\par{\xdef##1{\number\refnumber}%
\global\advance\refnumber by 1}%
\def\key##1##2\par{\expandafter\xdef%
\csname##1\endcsname{\number\refnumber}%
\global\advance\refnumber by 1}#1\par}}
\long\def\references{%
\penalty-800\vskip-\lastskip\vskip 15pt plus10pt minus5pt 
{\large\bf References}\ppar %`References' is set \large bold with space around.
{\leftskip=25pt\frenchspacing    % The list of references is set 
\small\parskip=3pt plus2pt       % \small  with small spaces between,
\def\refkey##1##2\par{\noindent  % numbers in [,]'s and set just to the
\llap{[##1]\stdspace}\ignorespaces##2\par}         % left of a 25pt margin.
\def\key##1##2\par{\noindent  
\llap{[\ref{##1}]\stdspace}\ignorespaces##2\par}  
\def\,{\thinspace}\thereflist\par}}
%
%   Next a footnote macro (with automatic numbering) which sets the
%   footnote  \small.
%
%   Typical use :
%         ..... are yellow.\fnote{By yellow here we mean Britsh
%    Standard colour BS3320.} 
%
\newcount\footnotenumber         % Register for footnote number
\footnotenumber=1                % set initially to 1
\def\fnote#1{\xdef\nextkey{\number\footnotenumber}%
{\small\ifnum\footnotenumber>9\parindent=14pt%
\else\parindent=10pt\fi\footnote{$^{\number\footnotenumber}$}%
{\hglue-5pt#1}\global\advance\footnotenumber by 1}}
%
%
%   Next macros for handling figures with automatic numbering (using 
%   TeX's \midinsert to float the figure to a suitable place).
%   
%   The \figure ... \endfigure macro centres the figure and adds
%   an automatically numbered label  Figure XX  after it.
%
%   If you have a caption, then type \caption{caption text} 
%   somewhere between \figure and \endfigure.  The macro
%   will then add  Figure XX: caption text  after the figure.
%
%   If you want an unnumbered or uncentred figure, then use TeX's raw 
%       \midinsert Figure instructions \endinsert  
%   and if you want a numbered figure label in the same style then
%   use \caption{caption text} outside of  \figure ... \endfigure.
%
%   If you need just the label Figure XX  outside of  \figure ... \endfigure
%   then type  \figurelabel .
%
\newcount\figurenumber          % register for figure number
\figurenumber=1                 % set initially to 1
\def\caption#1{\xdef\nextkey{\number\figurenumber}%
\cl{\small Figure \number\figurenumber: #1}%
\global\advance\figurenumber by 1}
\def\figurelabel{\xdef\nextkey{\number\figurenumber}%
\cl{\small Figure \number\figurenumber}%
\global\advance\figurenumber by 1}
\long\def\figure#1\endfigure{{\xdef\nextkey{\number\figurenumber}%
\let\captiontext\relax\def\caption##1{\xdef\captiontext{##1}}%
\midinsert\cl{\ignorespaces#1\unskip\unskip\unskip\unskip}\vglue6pt\cl{\small 
Figure \number\figurenumber\ifx\captiontext\relax\else: \captiontext
\fi}\endinsert\global\advance\figurenumber by 1}}
%
%   Macros for self-correcting internal references.
%
%   There are two macros  \key{KEY}  and  \ref{KEY} .
%
%   The \key macro sets up KEY as a key for whatever number is 
%   being referenced and the \ref macro converts the KEY into 
%   that number.  Type \key after a  \section or \proc or 
%   \label or \fnote or \figure or \caption or \figurelabel .
%
%   Example:
%
%       \section{Introduction}\key{intro}
%       \proc{Theorem}\key{MainTh}Lemons are yelloy\endproc
%       Here we follow\fnote{Follow in the sense of Dickens}
%       \key{Dickens-note}the crowd ....  
%
%       In section \ref{intro}
%       we stated theorem \ref{mainTh} and noted (see footnote 
%       \ref{Dickens-note}) ...
%
\def\nextkey{??}   %  initialise \nextkey (which is reset by all the
%                     numbering macros)
%
\def\key#1{\expandafter\xdef\csname #1\endcsname{\nextkey}}
\def\ref#1{\expandafter\ifx\csname #1\endcsname\relax
\immediate\write16{Reference {#1} undefined}??\else
\csname #1\endcsname\fi}
%
%   Note:  If the KEY contains only letters then \KEY has exactly the
%   same meaning as \ref{KEY} so in the example you could have:
%
%       In section \intro\ we ....
%
%   The \key will work at any time after the macro which sets the
%   number, provided no other macro which sets a number has been used. 
%
%   Macros for forward references:
%              =======
%   The \key \ref macros ONLY work for backwards references.  If you  
%   want to use forwards references, then type \useforwardrefs  near
%   the beginning of your file.  The KEY's are then stored in an
%   auxiliary  .ref  file and you then suffer the same disadvantage as
%   when using LaTeX that you must TeX the file twice to get
%   the references correct.
%
%   To use a forward ref type \ref{KEY}.  (You can type the
%   alternative  \KEY  but you'll get an error on first TeX'ing 
%   if the \KEY is not yet defined.) 
%
%   The macro also allows external references to be listed at the end 
%   of the file (if you wish to).  (Indeed they can be typed anywhere
%   before the \references command.)  You can combine the reference list
%   and the \references command by typing the references (using the
%   same syntax as before) between the commands \biblio and \endbiblio 
%   (don't type \references or they'll be printed twice).
%
\newread\gtinfile
\newwrite\gtreffile
\def\useforwardrefs{
\openin\gtinfile\jobname.ref
\ifeof\gtinfile
\closein\gtinfile
\immediate\write16{No file \jobname.ref}
\else
\closein\gtinfile
\input \jobname.ref
\fi
\immediate\openout\gtreffile \jobname.ref
%
%   Adapt \key :
%
\def\key##1{{\def\\{\noexpand}%
\expandafter\xdef\csname ##1\endcsname{\nextkey}%
\immediate\write\gtreffile{\\\expandafter\\\def\\\csname ##1\\\endcsname%
{\nextkey}}}}
%
%  Adapt macros for external references:  
%
\long\def\reflist##1\endreflist{%
\long\def\thereflist{##1}{\def\refkey####1####2\par{\xdef####1{%
\number\refnumber}{\def\\{\noexpand}\immediate\write\gtreffile
{\\\def\\####1{\number\refnumber}}}\global\advance\refnumber by 1}%
\def\key####1####2\par{\expandafter\xdef%
\csname####1\endcsname{\number\refnumber}%
{\def\\{\noexpand}\immediate\write\gtreffile
{\\\expandafter\\\def\\\csname ####1\\\endcsname{\number\refnumber}}}
\global\advance\refnumber by 1}##1\par}}
\long\def\biblio##1\endbiblio{\reflist##1\endreflist\references}%
%
%  Adapt obselete key macros (\numkey, \seckey and \figkey):
%
\def\numkey##1{{\def\\{\noexpand}%
\xdef##1{\number\sectionnumber.\number\resultnumber}
\immediate\write\gtreffile{\\\def\\##1%
{\number\sectionnumber.\number\resultnumber}}}}
\def\seckey##1{{\def\\{\noexpand}\xdef##1{\number\sectionnumber}
\immediate\write\gtreffile{\\\def\\##1{\number\sectionnumber}}}}
\def\figkey##1{\xdef##1{\number\figurenumber}%
{\def\\{\noexpand}\immediate\write\gtreffile%
{\\\def\\##1{\number\figurenumber}}}
\number\figurenumber\global\advance\figurenumber by 1}
}   %  end of \useforwardrefs
%
%
%   The next five macros are obselete and have been superseeded by
%   the general \key macro above.  They are included merely to 
%   maintain backward compatibility for the package:
%
%
\def\figkey#1{\xdef#1{\number\figurenumber}%
\number\figurenumber\global\advance\figurenumber by 1}
\def\fig#1#2\endfig{%
\midinsert\cl{#2}\vglue6pt\cl{\small Figure #1}\endinsert}
\def\newfig{\number\figurenumber\global\advance\figurenumber by 1}
\def\numkey#1{\xdef#1{\number\sectionnumber.\number\resultnumber}}
\def\seckey#1{\xdef#1{\number\sectionnumber}}
%
%   End of obselete macros.
%
%
%   The next macro is a version of the verbatim macro given by Knuth.
%
%   This macro produces a "verbatim" printout of
%   any ASCII string which does not contain the symbol "
%   (TeX files do not usually contain " 's).
%   More precisely, everything between consecutive pairs
%   of " 's is printed verbatim in the typewriter font cmtt.
%   For an explanation of how the macro works, see Knuth pp 420-1.
%
%   There are two switches: \verb (which switches the macro on)
%   and \brev which switches the macro off (the default).  When
%   the macro is switched off the symbol " has its usual 
%   meaning for TeX.  To use the macro, type \verb before use
%   and the use " to switch verbatim on and off.  Be careful
%   not to use " for any other purpose.  There is no need to
%   switch the macro off again unless you need to use " for
%   some other purpose (eg making  \mathchardef 's).  Note 
%   that the macro MUST BE OFF before inputting  amsnames.tex .
%
%   Whether the macro is on or off you can always use the
%   control sequence \dq (double quote) for " e.g.
%   \mathchardef\sum=\dq1350  is perfectly valid.
%   The control sequence \ttq is an abbreviation for
%   {\tt\dq}.  Thus "\ttq" will produce " (in cmtt)
%   inside a verbatim quote.
%
%
   %  define a code for " so it can be used when \verb is on
  %  code for " in cmtt
%
\def\verb{\catcode`\"=\active}       %  The main
\def\brev{\catcode`\"=12}            %  switches.
\brev                                %  Prime switches and
\verb                                %  switch on.
{\obeyspaces\gdef {\ }}              
{\catcode`\`=\active\gdef`{\relax\lq}}
\def"{%
\begingroup\baselineskip=12pt\def\par{\leavevmode\endgraf}%
\tt\obeylines\obeyspaces\parskip=0pt\parindent=0pt%
\catcode`\$=12\catcode`\&=12\catcode`\^=12\catcode`\#=12%
\catcode`\_=12\catcode`\~=12%
\catcode`\{=12\catcode`\}=12\catcode`\%=12\catcode`\\=12%
\catcode`\`=\active\let"\endgroup}
\brev      %   Finally switch the macro off (for safety)
%
%   Macros for itemised lists.   Typical use :
%    
%    \items
%    \item{(i)}Colours must be defined.
%    \item{(ii)}Colour cards may not be cited.
%    \enditems
%
%   Result :
%
%    (i)  Colours must be defined. 
%   (ii)  Colour cards may not be cited.
%
%
\def\items{\par\leftskip = 25pt}           % Start of itemised list         
\def\enditems{\par\leftskip = 0pt}         % end of itemised list   
\def\item#1{\par\leavevmode\llap{#1\stdspace}%
\ignorespaces}                             % labelled item
               % bulleted item.
%
%   The \quote ... \endquote macros are for typesetting quotations :
%

%
%   A few useful abbreviations :
%
\def\co{\colon\thinspace}    %  Colon with correct spacing for maps.
\def\np{\vfil\eject}         %  Forced page break (new page).
\def\nl{\hfil\break}         %  New line.
\def\cl{\centerline}         %  Centerline
\def\gt{{\mathsurround=0pt\it $\cal G\mskip-2mu$eometry \&\ 
$\cal T\!\!$opology}}        %  The journal title in recommended style
\def\gtm{{\mathsurround=0pt\it $\cal G\mskip-2mu$eometry \&\ 
$\cal T\!\!$opology $\cal M\mskip-1mu$onographs}}    %  for monographs
\def\agt{{\mathsurround=0pt\it$\cal A\mskip-.7mu$lgebraic \&\ 
$\cal G\mskip-2mu$eometric $\cal T\!\!$opology}}  % AGT
%
%    Finally some macros for automatic title page or header generation.
%    To use them type your header information using the following  
%    example as a guide :
%
%    Note that \\ is used as standard separator (for lines in \title and
%    \address, between authors and between email addresses or URL's)
%    and that \email, \url and \secondaddress are optional.
%

% Example:  \title{A short spoof paper\\with a two-line title}
% =======   \authors{Albert Einstein\\Leonardo da Vinci}
%           \address{IAS\\Princeton}\secondaddress{Renaissance\\Venice}
%           \email{ae@ias.princeton.edu\\ldv@ren.ven.hist}
%           \abstract 
%           A short spoof paper with a very short abstract.
%           \endabstract 
%           \primaryclass{00-01, 00-02}\secondaryclass{68-00, 68-01}
%           \keywords{Short, spoof, paper}
%           \maketitlepage
%
%
%    The title page or header will then be generated automatically.
%
%
%    Define the various ingredients of the title page:
%
\def\title#1{\def\thetitle{#1}}

\def\author#1{\edef\previousauthors{\theauthors}
 \ifx\theauthors\relax\def\theauthors{#1}\else
 \def\theauthors{\previousauthors\par#1}\fi}

\let\authors\author        % aliases
\def\address#1{\edef\previousaddresses{\theaddress}
 \ifx\theaddress\relax\def\theaddress{#1}\else
 \def\theaddress{\previousaddresses\par\vskip 2pt\par#1}\fi}
                             % alias
\def\secondaddress#1{\edef\previousaddresses{\theaddress}
 \ifx\theaddress\relax\def\theaddress{#1}\else
 \def\theaddress{\previousaddresses\par{\rm and}\par#1}\fi}   

\def\email#1{\edef\previousemails{\theemail}
 \ifx\theemail\relax\def\theemail{#1}\else
 \def\theemail{\previousemails\hskip 0.75em\relax#1}\fi}
  % aliases
\def\secondemail#1{\edef\previousemails{\theemail}
 \ifx\theemail\relax\def\theemail{#1}\else
 \def\theemail{\previousemails\hskip 0.75em{\rm and}\hskip 0.75em
 \relax#1}\fi}
\def\url#1{\edef\previousurls{\theurl}
 \ifx\theurl\relax\def\theurl{#1}\else
 \def\theurl{\previousurls\hskip 0.75em\relax#1}\fi}
      % aliases
\def\secondurl#1{\edef\previousurls{\theurl}
 \ifx\theurl\relax\def\theurl{#1}\else
 \def\theurl{\previousurls\hskip 0.75em{\rm and}\hskip 0.75em
 \relax#1}\fi}
\long\def\abstract#1\endabstract{\long\def\theabstract{#1}}
\def\primaryclass#1{\def\theprimaryclass{#1}}
                        % alias
\def\secondaryclass#1{\def\thesecondaryclass{#1}}
\def\keywords#1{\def\thekeywords{#1}}
%
%  Set \\ to \par and title page items to \relax to initialise macros :
%
\let\\\par\let\thetitle\relax\let\theshorttitle\relax
\let\theauthors\relax\let\theshortauthors\relax
\let\theaddress\relax\let\theshortaddress\relax
\let\theemail\relax\let\theurl\relax
\let\theabstract\relax\let\theprimaryclass\relax
\let\thesecondaryclass\relax\let\thekeywords\relax
%
%
%
%   Basic title page layout (edit this macro if you
%   wish to adjust the title page layout) :
%
\long\def\maketitlepage{    % start of definition of \maketitlepage

\vglue 0.2truein   % top margin

% title :
%
{\parskip=0pt\leftskip 0pt plus 1fil\def\\{\par\smallskip}{\large
\bf\thetitle}\par\medskip}   

\vglue 0.15truein 

% authors :
%
{\parskip=0pt\leftskip 0pt plus 1fil\def\\{\par}{\sc\theauthors}
\par\medskip}%
 
\vglue 0.1truein 

% address(es) email's and URL's (with switches to detect whether the
% optional items have been used) :
%
{\small\parskip=0pt
{\leftskip 0pt plus 1fil\def\\{\par}{\sl\theaddress}\par}
\ifx\theemail\relax\else  % email address?
\vglue 5pt \def\\{\stdspace{\rm and}\stdspace} 
\cl{Email:\stdspace\tt\theemail}\fi
\ifx\theurl\relax\else    % URL given?
\vglue 5pt \def\\{\stdspace{\rm and}\stdspace} 
\cl{URL:\stdspace\tt\theurl}\fi\par}

\vglue 7pt 

{\bf Abstract}

\vglue 5pt

\theabstract

\vglue 7pt 

{\bf AMS Classification numbers}\quad Primary:\quad \theprimaryclass\par

Secondary:\quad \thesecondaryclass

\vglue 5pt 

{\bf Keywords:}\quad \thekeywords

\np  % page break at the end of the title page

}    % end of definition of \maketitlepage
%
%    % \makeshorttitle (for general preprints) doesn't take a new page
%
\long\def\makeshorttitle{    % start of definition of \makeshorttitle

%\vglue 0.2truein   % top margin

% title :
%
{\parskip=0pt\leftskip 0pt plus 1fil\def\\{\par\smallskip}{\large
\bf\thetitle}\par\medskip}   

\vglue 0.05truein 

% authors :
%
{\parskip=0pt\leftskip 0pt plus 1fil\def\\{\par}{\sc\theauthors}
\par\medskip}%
 
\vglue 0.03truein 

% address(es) email's and URL's (with switches to detect whether the
% optional items have been used) :
%
{\small\parskip=0pt
{\leftskip 0pt plus 1fil\def\\{\par}{\sl\ifx\theshortaddress\relax
\theaddress\else\theshortaddress\fi}\par}
\ifx\theemail\relax\else  % email address?
\vglue 5pt \def\\{\stdspace{\rm and}\stdspace} 
\cl{Email:\stdspace\tt\theemail}\fi
\ifx\theurl\relax\else    % URL given?
\vglue 5pt \def\\{\stdspace{\rm and}\stdspace} 
\cl{URL:\stdspace\tt\theurl}\fi\par}

\vglue 10pt 

% abstract and classification numbers (with switches):

{\small\leftskip 25pt\rightskip 25pt{\bf Abstract}\stdspace\theabstract

{\bf AMS Classification}\stdspace\theprimaryclass
\ifx\thesecondaryclass\relax\else; \thesecondaryclass\fi\par
{\bf Keywords}\stdspace \thekeywords\par}
\vglue 7pt
}    % end of definition of \makeshorttitle
\let\maketitle\makeshorttitle        %% alias
%
%    %%%% \makeagttitle (for AGT) and \makemontitle
%         These are similar to \makeshorttitle but
%         with addresses omitted (they go at the end)
%
%%%% publication info and test defaults:

\def\volumenumber#1{\def\thevolumenumber{#1}}
\def\volumename#1{\def\thevolumename{#1}}
\def\volumeyear#1{\def\thevolumeyear{#1}}
\def\pagenumbers#1#2{\def\startpage{#1}\def\finishpage{#2}}
\def\published#1{\def\publishdate{#1}}

%% Defaults for authors to use to check layout
\volumenumber{X}
\volumename{Volume name goes here}
\volumeyear{20XX}
\pagenumbers{1}{XXX}
\published{XX Xxxember 20XX}

\long\def\makegtmontitle{   % start of definition of \makegtmontitle

\count0=\startpage

\gtm\nl        %   GT mongraphs (top left) 
{\small Volume \thevolumenumber: \thevolumename\nl 
Pages \startpage--\finishpage\nl}

\vglue 0.1truein   % top margin

% title
{\parskip=0pt\leftskip 0pt plus 1fil\def\\{\par\smallskip}{\large
\bf\thetitle}\par\medskip}   
\vglue 0.05truein 

% authors :
%
{\parskip=0pt\leftskip 0pt plus 1fil\def\\{\par}{\sc\theauthors}
\par\medskip}%
 
\vglue 0.03truein 

%  abstract and classification numbers:

{\small\leftskip 25pt\rightskip 25pt{\bf Abstract}\stdspace\theabstract

{\bf AMS Classification}\stdspace\theprimaryclass
\ifx\thesecondaryclass\relax\else; \thesecondaryclass\fi\par
{\bf Keywords}\stdspace \thekeywords\par}\vglue 7pt

}   % end of definition of \makegtmontitle

\long\def\makeagttitle{   %%% start of definition of \makeagttitle
\agt\hfill      %   Journal title (top left) 
%   logo placeholder (top right)
\hbox to 60truept{\vbox to 0pt{\vglue -14truept{\bf [Logo here]}\vss}\hss}
\break
{\small Volume \thevolumenumber\ (\thevolumeyear)
\startpage--\finishpage\nl
Published: \publishdate}

\vglue .2truein

% title
{\parskip=0pt\leftskip 0pt plus 1fil\def\\{\par\smallskip}{\large
\bf\thetitle}\par\medskip}   
\vglue 0.05truein 

% authors :
%
{\parskip=0pt\leftskip 0pt plus 1fil\def\\{\par}{\sc\theauthors}
\par\medskip}%
 
\vglue 0.03truein 

%  abstract and classification numbers:

{\small\leftskip 25truept\rightskip 25truept{\bf Abstract}\stdspace\theabstract

{\bf AMS Classification}\stdspace\theprimaryclass
\ifx\thesecondaryclass\relax\else; \thesecondaryclass\fi\par
{\bf Keywords}\stdspace \thekeywords\par}\vglue 7truept

}   %%%% end of definition of \makeagttitle

%%%%% Macro to typeset addresses (typically at the end of the paper)

\def\Addresses{\bigskip
{\small \parskip 0pt \leftskip 0pt \rightskip 0pt plus 1fil \def\\{\par}
\sl\theaddress\par\medskip \rm Email:\stdspace\tt\theemail\par
\ifx\theurl\relax\else\smallskip \rm URL:\stdspace\tt\theurl\par\fi}}

\def\agtart{%   Mock-up of AGT article style (for authors to test with)
%  get print centerpage:
\hoffset 14truemm
\voffset 31truemm
\font\phead=cmsl9 scaled 950
\font\pnum=cmbx10 scaled 913
\font\pfoot=cmsl9 scaled 950
%  headline and footline
\headline{\vbox to 0pt{\vskip -4.5mm\line{\small\phead\ifnum
\count0=\startpage ISSN numbers are printed here
\hfill {\pnum\folio}\else\ifodd\count0\def\\{ }% 
\ifx\theshorttitle\relax\thetitle\else\theshorttitle\fi\hfill{\pnum\folio}
\else\def\\{ and }{\pnum\folio}\hfill\ifx\theshortauthors\relax\theauthors
\else\theshortauthors\fi\fi\fi}\vss}}
\footline{\vbox to 0pt{\vglue 0mm\line{\small\pfoot\ifnum\count0=\startpage
Copyright declaration is printed here\hfill\else
\agt, Volume \thevolumenumber\ (\thevolumeyear)\hfill\fi}\vss}}
%  force \agttitle
\let\maketitle\makeagttitle\let\makeshorttitle\makeagttitle
\let\maketitlepage\makeagttitle}

\def\gtmonart{%   Mock-up of GT monograph style (for authors to test with)
%  get print centerpage:
\hoffset 14truemm
\voffset 31truemm
\font\phead=cmsl9 scaled 950
\font\pnum=cmbx10 scaled 913
\font\pfoot=cmsl9 scaled 950
%  headline and footline
\headline{\vbox to 0pt{\vskip -4.5mm\line{\small\phead\ifnum
\count0=\startpage ISSN numbers are printed here
\hfill {\pnum\folio}\else\ifodd\count0\def\\{ }% 
\ifx\theshorttitle\relax\thetitle\else\theshorttitle\fi\hfill{\pnum\folio}
\else\def\\{ and }{\pnum\folio}\hfill\ifx\theshortauthors\relax\theauthors
\else\theshortauthors\fi\fi\fi}\vss}}
\footline{\vbox to 0pt{\vglue 0mm\line{\small\pfoot\ifnum\count0=\startpage
Copyright declaration is printed here\hfill\else
\gtm, Volume \thevolumenumber\ (\thevolumeyear)\hfill\fi}\vss}}
%  force \makegtmontitle
\let\maketitle\makegtmontitle\let\makeshorttitle\makegtmontitle
\let\maketitlepage\makegtmontitle}

\def\gtart{%   Mock-up of GT article style (for authors to test with)
%  get print centerpage:
\hoffset 14truemm
\voffset 31truemm
\font\phead=cmsl9 scaled 950
\font\pnum=cmbx10 scaled 913
\font\pfoot=cmsl9 scaled 950
%  headline and footline
\headline{\vbox to 0pt{\vskip -4.5mm\line{\small\phead\ifnum
\count0=\startpage ISSN numbers are printed here
\hfill {\pnum\folio}\else\ifodd\count0\def\\{ }% 
\ifx\theshorttitle\relax\thetitle\else\theshorttitle\fi\hfill{\pnum\folio}
\else\def\\{ and }{\pnum\folio}\hfill\ifx\theshortauthors\relax\theauthors
\else\theshortauthors\fi\fi\fi}\vss}}
\footline{\vbox to 0pt{\vglue 0mm\line{\small\pfoot\ifnum\count0=\startpage
Copyright declaration is printed here\hfill\else
\gt, Volume \thevolumenumber\ (\thevolumeyear)\hfill\fi}\vss}}
%  force \maketitlepage
\let\maketitle\maketitlepage\let\makeshorttitle\maketitlepage}

\input epsf
%
% Initialise for MTG control
%
%
\def\relabelbox{%
  \hbox\bgroup%
%
%
}%
\def\endrelabelbox{%
\egroup%
}%
\def\relabel #1#2 {%
  \special{ps:/a {} def}%
%
  \smash{\rlap{#2}}%
%
}%
\def\adjustrelabel <#1,#2> #3#4 {%
  \special{ps:/a {} def}%
%
  \smash{\rlap{\kern #1 \raise #2\hbox{#4}}}%
%
}%
\def\extralabel <#1,#2> #3 {\smash{\rlap{\kern #1 \raise #2\hbox{#3}}}}%

\input xyall
\hoffset 14truemm
\voffset 31truemm

\def\strut{\vrule width 0pt height 11pt}

%%% abbreviations etc

\def\gen#1{\langle #1 \rangle}
\def\ngen#1{\langle \langle #1 \rangle \rangle}
\let\la\langle
\let\ra\rangle

\let\bar\overline
\font\spec=cmtex10 scaled 1095 
\def\d{\hbox{\spec \char'017\kern 0.05em}}

\def\ex{\mathop{\rm ex}}
\def\inv{^{-1}}
\def\co{\mskip 1.5mu\colon\thinspace}

\def\o{{\cal O}}
\def\Z{\Bbb Z}
\let\wtilde\widetilde
\let\what\widehat

%%%%%%%%%%%%
% get figure captions to wrap correctly:

\long\def\figuresc#1\endfigure{{\xdef\nextkey{\number\figurenumber}%
\let\captiontext\relax\def\caption##1{\xdef\captiontext{##1}}%
\midinsert\cl{\ignorespaces#1\unskip\unskip\unskip\unskip}\vglue6pt{\small 
\parskip=0pt\leftskip 0pt plus 1fil\def\\{\par}Figure 
\number\figurenumber\ifx\captiontext\relax\else: \captiontext
\fi\par}\endinsert\global\advance\figurenumber by 1}}

\long\def\figure#1\endfigure{{\xdef\nextkey{\number\figurenumber}%
\let\captiontext\relax\def\caption##1{\xdef\captiontext{##1}}%
\midinsert\cl{\ignorespaces#1\unskip\unskip\unskip\unskip}\vglue6pt{\small 
\parskip=0pt\leftskip 25pt \rightskip 25pt\def\\{\par}Figure 
\number\figurenumber\ifx\captiontext\relax\else: \captiontext
\fi\par}\endinsert\global\advance\figurenumber by 1}}

\reflist

\key{BBP}
{\bf Y G Baik}, {\bf W A Bogley}, {\bf S J Pride}, {\it On the
asphericity of length four relative group presentations},
Internat. J. Algebra Comput. 7 (1997) 277--312 

\key{BoPr}
{\bf W A Bogley}, {\bf S J Pride}, {\it Aspherical relative
presentations}, Proc. Edinburgh Math. Soc. 35 (1992) 1--39 

\key{Brown}
{\bf Kenneth S Brown}, {\it Cohomology of Groups}, Springer-Verlag (1982) 

\key{CR}
{\bf Marshall M Cohen}, {\bf Colin Rourke}, {\it The surjectivity problem
for one-gener\-at\-or, one-relator extensions of torsion-free groups},
Geom. Topol. 5 (2002) 127--142

\key{Edj}
{\bf Martin Edjvet}, {\it On the asphericity of one-relator relative
presentations}, Proc. Roy. Soc. Edinburgh Sect. A 124 (1994) 713--728 

\key{FR}
{\bf Roger Fenn}, {\bf Colin Rourke}, {\it Klyachko's methods and the
solution of equations over torsion-free groups}, L'Enseignment
Math\'ematique, 42 (1996) 49--74

\key{FRo}
{\bf Roger Fenn}, {\bf Colin Rourke}, {\it Characterisation of a class
of equations with solutions over torsion-free groups}, from: ``The
Epstein Birthday Schrift'', (I Rivin, C Rourke and C Series, editors),
Geometry and Topology Monographs, Volume 1 (1998) 159-166

\key{Ho}
{\bf James Howie}, {\it The solution of length three equations over
groups}, Proc. Edinburgh Math. Soc. 26 (1983) 89--96

\key{HoMet}
{\bf J Howie}, {\bf V Metaftsis}, {\it On the asphericity of length five
relative group presentations}, Proc. London Math. Soc. (3) 82 (2001)
173--194 

\key{Kl}
{\bf Anton A Klyachko}, {\it A funny property of sphere and equations over
groups}, Comm. in Alg. 21 (1993) 2555--2575

\key{Lev}{\bf Frank Levin}, {\it Solutions of equations over groups},
Bull. Amer. Math. Soc. 68 (1962) 603--604 

\key{Met}
{\bf V Metaftsis}, {\it On the asphericity of relative group
presentations of arbitrary length}, Internat. J. Algebra Comput. 13 no. 3
(2003) 323--339 

\key{Dunce}
{\bf Colin Rourke}, {\it On dunce hats and the Kervaire conjecture},
Papers presented to Christopher Zeeman, University of Warwick (1988)
221--230, available from: {\tt
http://www.maths.warwick.ac.uk/\char'176cpr/ftp/dunce.ps}

\key{Short} {\bf H Short}, {\it	 Topological methods in group theory: the
adjunction problem},  Warwick Ph.D. thesis (1983)

\key{Whitehead} {\bf J\,H\,C Whitehead}, {\it On adding relations to homotopy
groups}, Ann. of Math. 42 (1941) 409--428

\endreflist

\title{Diagrams and the second homotopy group}
\authors{Max Forester\fnote{Research supported by EPSRC grant 
GR/N20867}\\Colin Rourke}
\address{Mathematics Institute, University of Warwick,
Coventry, CV4 7AL, UK}
\email{forester@maths.warwick.ac.uk, cpr@maths.warwick.ac.uk}

\abstract
We use Klyachko's methods [\CR,\FR,\FRo,\Kl] to prove that, if a
1--cell and a 2--cell are added to a complex with torsion-free
fundamental group, and with the 2--cell attached by an amenable
$t$--shape, then $\pi_2$ changes by extension of scalars.  It then
follows using a result of [\BoPr] that the resulting fundamental group
is also torsion free. We also prove that the normal closure of the
attaching word contains no words of smaller complexity.
\endabstract

\primaryclass{57M20, 57Q05}\secondaryclass{20E22, 20F05}
\keywords{Diagram, second homotopy group, Whitehead conjecture, 
torsion-free group, amenable shape, Kervaire conjecture}

\maketitle

\medskip

\rk{\large\bf Introduction}

\medskip

We consider the following problem.  Suppose that $K$ is a CW complex
and that $L$ is formed from $K$ by adding cells of dimension $\le2$.

\proclaim{Problem}\rm  Describe $\pi_2(L)$ in terms of $\pi_2(K)$.
\endproc 

The problem is interesting in its own right, but also because of 
its relation to Whitehead's conjecture [\Whitehead] (still
unsolved).  This says that if $K$ is 2--dimensional with non-zero
$\pi_2$ then $L$ also has non-zero $\pi_2$.

We shall answer the question completely in a special case.  Suppose
that $L=K\cup e^1\cup e^2$ and let $t$ be the new generator of $\pi_1$
determined by $e^1$.  Then the attaching map for $e^2$ represents a
word $w$ in $G*\gen{t}$ where $G=\pi_1(K)$ which we can assume to be
cyclically reduced. The {\sl $t$--shape} of $w$ is the unreduced word
formed by occurrences of $t$. There is a useful notion of {\sl
amenable $t$--shape} that was introduced in [\FR]. This includes all
$t$--shapes having total exponent $\pm 1$ in $t$. The exact definition
will be recalled in section 1.

\proclaim{Main Theorem} 
Suppose that $L=K\cup e^1\cup e^2$ and that the attaching map for
$e^2$ represents $w$ in $G*\gen{t}$ where $G=\pi_1(K)$ and $t$ is the
new generator of $\pi_1$ determined by $e^1$.  Suppose\eject
\items
\item{\rm(1)}$G$ is
torsion-free 
\item{\rm(2)}the $t$--shape of $w$ is amenable. 
\enditems
Then there is an isomorphism of $\Z
\pi_1(L)$--modules $$\pi_2 (L) \ \cong \ \Z \pi_1 (L) \otimes_{\Z \pi_1
(K)} \pi_2 (K).$$ 
\endproc

This shows that in passing from $K$ to $L$, $\pi_2$ changes in the
simplest possible way, ie by ``extension of scalars''. No elements are
killed, and all new elements are accounted for by the change in
fundmental group. 

Note that under conditions (1) and (2), Fenn and Rourke have solved
the adjunction problem, and proved that $\pi_1(K)\to\pi_1(L)$ is
injective [\FR].  From this it follows easily that if $K$ is
2--dimensional then $\pi_2(K)\to\pi_2(L)$ is also injective.  (Let
$\wtilde{K}, \wtilde{L}$ be the universal covers; then $K\subset L$
lifts to $\wtilde{K}\subset\wtilde{L}$ and hence $\pi_2(K)\cong
H_2(\wtilde{K})\subset H_2(\wtilde{L})\cong \pi_2(L)$, where the 
middle inclusion follows from the fact that $\wtilde L-\wtilde K$
is 2--dimensional.)  So this paper gives no new information for
the Whitehead conjecture itself.  However note that we do not need to
assume that $K$ is 2--dimensional for our results.

The key observation of this paper is that the Klyachko proof used by
Fenn and Rourke actually shows much more than is needed for the
adjuction problem.  For the latter one needs that there are no
diagrams based on $w$ with a non-trivial boundary.  The proof actually
shows that {\sl there are no irreducible diagrams\/} whatever the
boundary.  This translates into the statement that $\pi_2(K^+,K)$ maps
onto $\pi_2(L,K)$ where $K^+=K\cup e^1$.  From this the Main Theorem
follows by standard algebraic topology.

It is interesting to note that, whilst the adjunction problem is still
open for torsion-free groups, and indeed for general groups provided
the exponent sum of $t$ in $w$ is non-zero, the Main Theorem fails if
either hypothesis (1) or (2) is dropped.  We shall see that in either
case there are both irreducible diagrams and elements of $\pi_2(L,K)$
not accessible from $\pi_2(K^+,K)$.  This has consequences for a
possible proof of the adjunction problem.  Fenn and Rourke [\FR; page
70] suggest that the Klyachko methods should be strong enough to prove
the 1--generator 1--relator adjunction problem for torsion-free groups
for any $t$--shape.  While this may be true, this paper shows that the
methods will need to be extended considerably if they are to work for
a $t$--shape which is a proper power.

This paper provides an alternate approach to part of the Klyachko
method, namely the ``algebraic trick'' described in [\FR; pages
64--66]. Here this step is explained in terms of diagrams; see in
particular figure 6. 

The diagram methods also apply to the Cohen--Rourke results [\CR] and
imply strong information about the normal closure $\ngen w$ of $w$ in
$G*\gen t\,$: no word of complexity smaller than the complexity of $w$
lies in $\ngen w$ (complexity is defined in section 1).  It should be
noted that this last result can be deduced directly from the
Cohen--Rourke methods without using the methods of the present paper.

The basic geometric result, that there are no irreducible diagrams, is
related to the notion of asphericity of relative presentations introduced
in [\BoPr] and also studied in [\Edj, \BBP, \HoMet, \Met]. However we
use a slightly different notion of irreducibility, involving
basepoints. This is to ensure that ``dipoles'' can always be removed by a
homotopy. The difference between the two notions is readily apparent in
the second example of section 4. 

On the other hand, the two notions differ only in the case of words that
admit a non-trivial symmetry (considered as cyclic words), and such words
are outside the scope of the Main Theorem. Hence the basic diagram result
may be interpreted as saying that relative group presentations satisfying
conditions (1) and (2) are aspherical in the sense of [\BoPr]. 

Combining our result with one of the main results of [\BoPr] it
follows that $\pi_1(L)$ is also torsion free, which answers a question of
Cohen and Rourke [\CR].  In a future paper we shall give a direct
proof of this result and deduce consequences for the multivariable
adjunction problem. We would like to thank the referee for pointing out
the connection of our result with the notion of aspherical relative
presentations, which has made this important consequence clear to us. 

Here is an outline of the paper.  In section 1 we define diagrams and
irreducibility and state our main results about the non-existence of
diagrams.  Section 2 contains the proofs of the diagram theorems and
in section 3 we translate the diagram results into statements about
$\pi_2$ and prove the Main Theorem. In section 4 we give the
counterexamples mentioned above and discuss the limits of the Klyachko
methods for proving the adjunction problem over torsion-free groups. 
In section 5 we give a generalisation of the Diagram Theorem which
implies the complexity result mentioned above.

\section{Diagrams}

Let $G$ be a group and $w\in G*\gen{t}$.  By a disc $D$ {\sl labelled
by $w$} we mean a 2--disc with legs, one for each occurrence of $t$ or
$t^{-1}$ in $w$.  The legs are ``thick arcs'', in other words trivial
bundles with fibre $I$.  Further the fibres are oriented and labelled
$t$ in such a way that the $t$--shape of $w$ can be read from the legs
by reading around the the boundary in an anti-clockwise direction and
reading $t$ if we cross a leg whose fibres are oriented in this
direction and $t\inv$ otherwise. Further the segments of the boundary
of $D$ between the legs are labelled by elements of $G$ and the whole
of $w$ is recovered by reading these labels and the legs. Finally we
fix a basepoint in the word $w$ (such as at the beginning) and
indicate this on the boundary of the $w$--disc.

A disc labelled by $\bar w$ means the reflection of a disc labelled by
$w$.  A {\sl $w$--diagram} means a finite collection of discs in the
plane labelled by $w$ and $\bar w$ together with thick arcs with
oriented fibres labelled $t$ (called {\sl $t$--arcs}) which complete
the legs compatibly with the orientations. Fibred annuli ({\sl
$t$--circles}) are also allowed. The diagram cuts the plane
into a number of finite regions (the ``inside regions'') and one
infinite region the ``outside''. 

Each region of the diagram determines a word in $G$ by reading
anticlockwise around the region for inside regions and clockwise for
the outside region using the labels on the adjacent discs.  We require
that this word is the identity in $G$ for each inside region.  An
example is shown in figure 1.  In this example $w=tatbt\inv c$ and
$G$ satisfies the relations $b^2=c^2=1$, $a^{-1}ba=c$. 

\figuresc
\relabelbox
\epsfbox{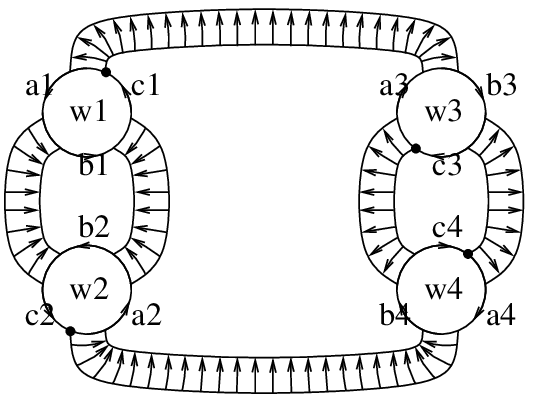}
\adjustrelabel <0mm, .2mm> {a1}{$\scriptstyle{a}$}
\adjustrelabel <-.1mm, .5mm> {a2}{$\scriptstyle{a}$}
\adjustrelabel <-.2mm, .2mm> {a3}{$\scriptstyle{a}$}
\adjustrelabel <-.3mm, .2mm> {a4}{$\scriptstyle{a}$}
\adjustrelabel <0mm, -.6mm> {b1}{$\scriptstyle{b}$}
\adjustrelabel <0mm, -.1mm> {b2}{$\scriptstyle{b}$}
\adjustrelabel <0mm, 0mm> {b3}{$\scriptstyle{b}$}
\adjustrelabel <.2mm, .2mm> {b4}{$\scriptstyle{b}$}
\adjustrelabel <0mm, 0mm> {c1}{$\scriptstyle{c}$}
\adjustrelabel <.2mm, 0mm> {c2}{$\scriptstyle{c}$}
\adjustrelabel <0mm, -.1mm> {c3}{$\scriptstyle{c}$}
\adjustrelabel <0mm, 0mm> {c4}{$\scriptstyle{c}$}
\adjustrelabel <0mm, 0mm> {w1}{$w$}
\adjustrelabel <0mm, 0mm> {w2}{$w$}
\adjustrelabel <0mm, -.2mm> {w3}{$\bar w$}
\adjustrelabel <0mm, -.2mm> {w4}{$\bar w$}
\endrelabelbox
\caption{A $w$--diagram with $w = tatbt^{-1}c$}
\endfigure

A diagram is {\sl irreducible} if:

(1)\qua it is connected,

(2)\qua it contains at least one disc labelled by $w$ or $\bar w$,

(3)\qua it does not contain pairs of discs of the type illustrated in
figure 2, where the joining `leg' represents the same occurrence of $t$
or $t\inv$ in $w$ and $\bar w$.

\figuresc
\relabelbox
\epsfbox{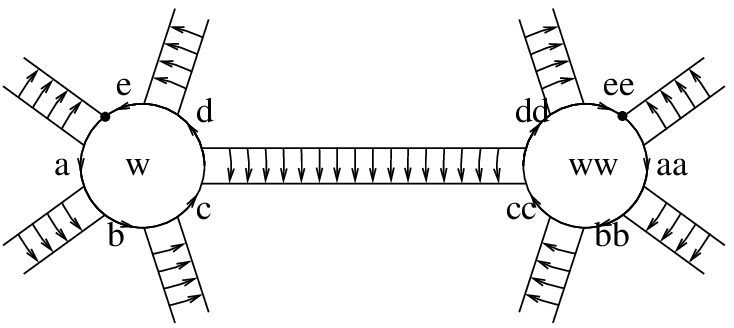}
\adjustrelabel <0mm, 0mm> {w}{$w$}
\adjustrelabel <0mm, -.2mm> {ww}{$\bar w$}
\adjustrelabel <.3mm, .3mm> {a}{$\scriptstyle{a}$}
\adjustrelabel <-.3mm, .3mm> {aa}{$\scriptstyle{a}$}
\adjustrelabel <.4mm, -.2mm> {b}{$\scriptstyle{b}$}
\adjustrelabel <1mm, -.2mm> {bb}{$\scriptstyle{b}$}
\adjustrelabel <0mm, .3mm> {c}{$\scriptstyle{c}$}
\adjustrelabel <.7mm, .3mm> {cc}{$\scriptstyle{c}$}
\adjustrelabel <0mm, -.5mm> {d}{$\scriptstyle{d}$}
\adjustrelabel <-.3mm, -.5mm> {dd}{$\scriptstyle{d}$}
\adjustrelabel <-.5mm, -.3mm> {e}{$\scriptstyle{e}$}
\adjustrelabel <-.1mm, -.4mm> {ee}{$\scriptstyle{e}$}
\endrelabelbox
\caption{Disc pair with $w = t^{-1}atbtct^{-1}dte$}
\endfigure

Condition (2) means that the diagram does not just consist of one
$t$--circle.  Condition (3) needs care if $w$ has symmetry.  The ``same
occurrence of $t$ or $t\inv$'' means the same position with respect to the
basepoint. 
 
The diagram shown in figure 1 is irreducible. 

We now recall the definition of
{\sl amenable\/} $t$--shapes from [\FR, \FRo].  Note that $t$--shapes
are always considered as cyclic words.

Write a given $t$--shape $T$ as an unreduced cyclic word in the symbols
$t$ and $t\inv$.  The {\sl Magnus derivative\/} $D(T)$ is the
$t$--shape given by deleting all subwords $tt^{-1}$ and closing up.  A
more symmetrical definition is given in [\FRo; page 161].  It is easy
to see that $D^m(T) = t^q$ for some $q\in \Z$ for $m$ sufficiently large.
Call min$\{m\mid D^m(T) = t^q\}$ the {\sl complexity\/} of $T$
and call $D^{m-1}(T)$ the {\sl root\/} of $T$.  A {\sl clump} in a
$t$--shape is a maximal connected subsequence of the form $t^q$ or
$t^{-q}$ where $q>1$.  A {\sl one-clump\/} shape is a $t$--shape which
has exactly one clump which is not the whole shape.  In other words it
is a shape of the form $t^p(tt^{-1})^q$ for $p,q >0$ or its inverse.
A $t$--shape is {\sl amenable\/} if it has root either $tt^{-1}$ or a
one-clump shape.  For future reference, we shall call these shapes
{\sl amenable root shapes}.

It is easy to see that any $t$--shape of total exponent sum 1 in $t$
is amenable since the root must be of the form $t(tt^{-1})^q$ for
$q>0$. On the other hand, no $t$--shape which is a proper power is
amenable. 

\proclaim{Diagram Theorem}
Suppose that $w$ is a word in $G*\gen{t}$ with an amenable $t$--shape and
that $G$ is a torsion-free group. Then there are no irreducible
$w$--diagrams. 
\endproc 

The proof is really an observation.  This is what the Klyachko proof
of the Kervaire conjecture (and the extension given by Fenn--Rourke in
[\FR, \FRo]) actually shows.  We shall need to go carefully through
the details to make this clear and we shall do this in the next
section.

If $G$ is not torsion-free, then there are many easily constructed
irreducible diagrams.  In [\Dunce] it is shown how to construct
infinite families of irreducible diagrams made of ``units''.  Figure 1
is a simple example of this construction: it is made of two 2--units.
All diagrams made in this way require the group to have torsion.
Short [\Short] has investigated diagrams made of units in detail and
proved that none gives a counterexample to the Kervaire conjecture.

If $G$ is torsion-free, but the word is not amenable, then again
irreducible diagrams may exist.  Explicit examples are given in section
4.

\section{Proof of the Diagram Theorem}
 
We follow the general pattern of proof of [\Kl] as exposited in [\CR,
\FR]. We prove first a special case, when the $t$--shape is an 
amenable root shape (either $tt\inv$ or a one-clump shape) and then
use the ``algebraic trick'' described in [\FR; pages 64--66] to
convert this to the general case.

We need to generalise diagrams for a set of words in $\Gamma \ast \langle
t \rangle$.  The Root Shape Theorem below yields the special case of the
Diagram Theorem by taking $\Gamma = G$. In the proof of the general case 
$\Gamma$ will represent $G \ast \langle s \rangle$. 

Let $\Gamma$ be a group and $w_0\in\Gamma*\gen t$ be cyclically reduced.
An element of $\Gamma$ in $w_0$ sandwiched cyclically between an
occurrence of $t$ and $t \inv$ is called a {\sl top coefficient} of
$w_0$. One sandwiched between an occurrence of $t\inv$ and $t$ is called
a {\sl bottom coefficient}.  The others are called {\sl middle
coefficients}.  

Now let $H$ be a subgroup of $\Gamma$ and $g$ an element of
$\Gamma$. We say that $g$ is {\sl free relative to $H$} if the subgroup
$\langle g, H\rangle$ of $\Gamma$ generated by $g$ and $H$ is naturally
the free product $\langle g \rangle \ast H$ of an infinite cyclic group
$\langle g \rangle$ with $H$. 

We shall use the following working hypotheses.  

\proclaim{Working hypotheses}

Suppose that $H$ and $H'$ are two isomorphic subgroups of a group
$\Gamma$ under the isomorphism $h\to h^{\phi}$, $h\in H$.  

Let $w_0$ be a word in $\Gamma*\gen t$ with amenable root $t$--shape
and let $a_i, b_i\in\Gamma$ be the bottom and top coefficients of $w_0$
repectively listed in any order.  Suppose that for each $i$, $a_i$ is
free relative to $H$ and $b_i$ is free relative to $H'$.  Let
$W\subset\Gamma*\gen t$ be the set of words $\{w_0,\,
h^t(h^{\phi})\inv\mid h\in H\}$ where $h^t=t\inv ht$. \endproc

\proclaim{Definitions}\rm  A {\sl $W\!$--diagram\/} is a diagram with discs
labelled by elements of $W$ or their inverses (denoted $\bar w_0$
etc) with the same interpretation as in section 1, together with arcs
completing the legs and such that words read around the boundaries of
inside regions are the identity in $\Gamma$, as in the definition of
a $w$--diagram.

Notice that a disc labelled by $h^t(h^{\phi})\inv$ has just two legs
emerging and the boundary arcs between the legs are labelled $h$ and
$h^\phi$; see figure 3.

A $W\!$--diagram is {\sl irreducible} if:

(1)\qua it is connected,

(2)\qua it contains at least one disc labelled by $w_0$ or $\bar w_0 $,

(3)\qua it does not contain pairs of discs of the type illustrated in
figure 2, where the joining `leg' represents the same occurrence of $t$
or $t\inv$ in $w_0$ and $\bar w_0 $,

(4)\qua it does not contain a string of two-leg discs labelled in order
by $h_i^t(h_i^{\phi})\inv, i=1,2,\ldots,p$ where $h_1h_2\ldots h_p=1$; 
see figure 3.

\figuresc
\relabelbox
\epsfbox{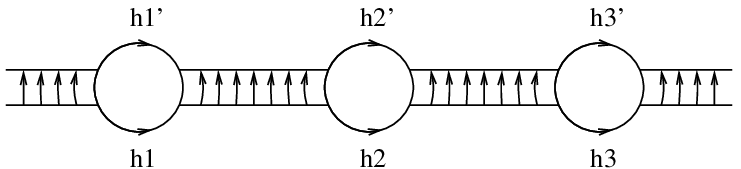}
\adjustrelabel <0mm, .5mm> {h1}{$\scriptstyle{h_1}$}
\adjustrelabel <0mm, .5mm> {h2}{$\scriptstyle{h_2}$}
\adjustrelabel <0mm, .5mm> {h3}{$\scriptstyle{h_3}$}
\adjustrelabel <0mm, 0mm> {h1'}{$\scriptstyle{h_1^{\phi}}$}
\adjustrelabel <0mm, 0mm> {h2'}{$\scriptstyle{h_2^{\phi}}$}
\adjustrelabel <0mm, 0mm> {h3'}{$\scriptstyle{h_3^{\phi}}$}
\endrelabelbox
\caption{A string of two-leg discs}
\endfigure

\proclaim{Root Shape Theorem}
Assume the working hypotheses.  Then there are no irreducible
$W\!$--diagrams.\endproc

Note that by taking $H$ and $H'$ to be trivial in the theorem we can
deduce a special case of the Diagram Theorem, namely the case when the
$t$--shape is an amenable root shape. 

\prf Let $D$ be an irreducible $W\!$--diagram.  We convert $D$ into 
a cell subdivision of the 2--sphere, which is essentially the dual
diagram, by putting a vertex in each region and joining by an edge
across each $t$--arc; see [\CR; figure 5].  We now define a traffic
flow as described on page 68 of [\FR] and obtain a contradiction in
exactly the same way as in [\CR] or [\FR].  The flow is described for
suitable $t$--shapes which are more general than one-clump shapes but
do not include the shape $tt\inv$, but this shape has a very easy
flow.  All cells are bigons and a traffic flow (without stops) is
defined as for the bigons on the right of [\FR; figure 5]. The details
for the contradiction are clearest in [\CR; page 137].  Crashes can
only occur at vertices.  There may be a crash at the vertex in the
outside region of $D$.  But since there are at least two crashes in
different places, there must be a crash at an internal vertex.  But
the flow has been chosen so that, at an internal vertex where all the cars
come together at the same time, the labels around the corners are all
$\{a,\, a\inv\}$ for some coefficient $a=a_i$ of $w_0$ together with
elements of $H$ or $\{b,\, b\inv\}$ for some coefficient $b=b_i$ of
$w_0$ together with elements of $H'$.  For definiteness assume that we
are in the former situation.  Then we can read an (unreduced) word of
the form $a^{\epsilon_1} h_1 h_2 \dots h_{i_1} a^{\epsilon_2} h_1 h_2
\dots h_{i_2}a^{\epsilon_3}\dots$ which is 1 in $\Gamma$.  Now if this
word contains a subword of the form $a^\epsilon a^{-\epsilon}$ then $D$
fails part (3) of the definition of irreducibility and if it contains a
subword of the form $h_1 h_2 \dots h_i$ which is 1 in $\Gamma$ then it
fails part (4).  Since $D$ is irreducible neither of these happen and the
word gives a non-trivial relation in $\gen{a,H}$, contradicting the
assumption that $a$ is free relative to $H$. 
\endprf

\rk{Proof of the Diagram Theorem}

We now convert the Root Shape Theorem into the Diagram Theorem by
using the ``algebraic trick'' of Klyachko described in [\FR; pages
64--66] and [\CR; pages 137--139].  However we need to do this on the
level of diagrams. Recall the following definitions:

Consider the exponent sum homomorphism $\ex\co G*\gen s \to \Z$.  It is
well known 
that $K$, the kernel of $\ex$, is a free product of copies of $G$
generated by elements of the form $g^{s^\o}=s^{-\o}gs^\o$, $1 \neq g\in
G$.

Any element of $K$ has a canonical expression of the form
$k=g_1^{s^{\o_1}}\cdots g_r^{s^{\o_r}}$, where $\o_i\ne\o_{i+1}$ for
each $i$.   We shall call the $g_i^{s^{\o_i}}$ the {\sl canonical
elements} of $k$.  Let min$(k)$ be the minimum value of $\o_i$,
$i=1,\ldots,r$ and  max$(k)$ the maximum value. Fix a positive integer
$m$. Consider the following subgroups of $K$:
$$\eqalign{%
H&=\la k\in K\mid \,\hbox{\rm min}(k)\ge0, \hbox{\rm max}(k)\le m-2\ra  \cr
H'&=\la k\in K\mid \,\hbox{\rm min}(k)\ge1, \hbox{\rm max}(k)\le m-1\ra \cr
J&=\la k\in K\mid \,\hbox{\rm min}(k)\ge0, \hbox{\rm max}(k)\le m-1\ra  \cr}$$
and the following subsets:
$$\eqalign{%
X&=\{ k\in K\mid \,\hbox{\rm min}(k)=0, \hbox{\rm max}(k)\le m-1\}  \cr
Y&=\{ k\in K\mid \,\hbox{\rm min}(k)\ge0, \hbox{\rm max}(k)=m-1\}  . \cr}$$

We consider words in $(G \ast \gen s) \ast \gen t$ obtained from an
amenable root $t$--shape by 
``blowing up''.  To be precise, we cyclically insert an element of $X$ (respectively $Y$, $J$) in each {\sl top} (respectively
{\sl bottom}, {\sl middle}) position, ie between successive occurrences
of $t$ and $t\inv$ (respectively $t\inv$ and $t$, $t$ and $t$ or $t\inv$
and $t\inv$). 

Now take the resulting word $w(s,t)$, substitute $t$ for $s$, and reduce
cyclically to obtain $w$. The $t$--shape of any $w$ obtained in this way
is {\sl amenable of complexity $m$}. See [\FRo] for the details of the
proof that this coincides with the definition given earlier and the
(apparently more general) one given in [\FR]. The root shape of $w$ is
the original $t$--shape before blowing up. We call $w(s,t)$ the {\sl
blown up form} of $w$ (also called normal form in [\FR,~\FRo]). Write it
as $$w(s,t)=t^{\epsilon_1} x_1(s)t^{\epsilon_2} x_2(s) \ldots
t^{\epsilon_n} x_n(s)$$ where each $x_i(s)$ belongs to the appropriate
subset $X$, $Y$ or $J$. 

We now describe a procedure for converting irreducible $w$--diagrams to
irreducible $W\!$--diagrams. Here $W = \{w(s,t), \, h^t(h^{\phi})\inv \mid
h \in H \}$ where $\phi \co H \to H'$ is given by $h^{\phi} = h^s$, and
$w(s,t)$ is regarded as a word with coefficients in $\Gamma = G \ast \gen
s$ and having amenable root $t$--shape. The process is illustrated in
figures 4--6 for the amenable word $w = attbt^{-1}t^{-1}cdt^{-1}et^2$. 

Given a $w$--diagram, first we replace every $w$--disc by a disc labelled
$w(s,t)$ together with additional $t$--arcs effecting the reduction from
$w(t,t)$ to $w$; see figure 4. The arcs in the outer ring are $t$--arcs,
but the legs inside are designated as $t$--legs or $s$--legs according to
$w(s,t)$. 

\figure
\relabelbox
\epsfbox{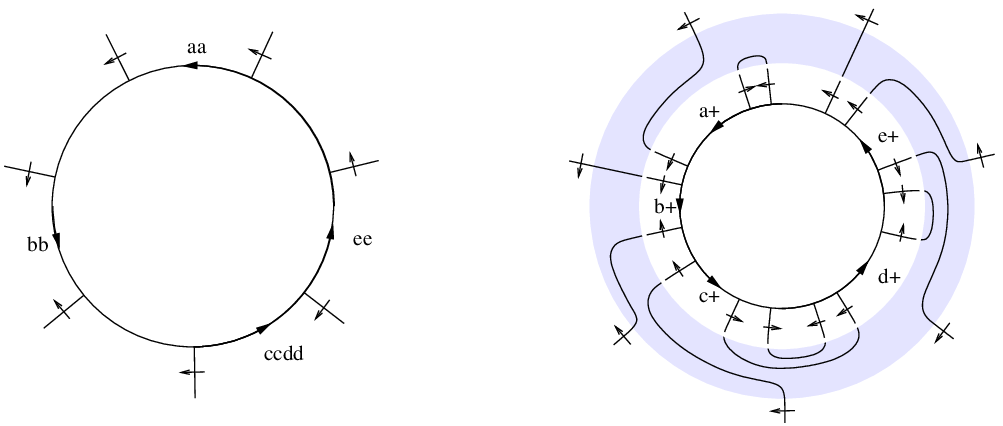}
\adjustrelabel <-.5mm, 0mm> {aa}{$\scriptstyle{a}$}
\adjustrelabel <0mm, 0mm> {bb}{$\scriptstyle{b}$}
\adjustrelabel <-.5mm, .5mm> {ccdd}{$\scriptstyle{cd}$}
\adjustrelabel <-1mm, .3mm> {ee}{$\scriptstyle{e}$}
\adjustrelabel <-.4mm, 0mm> {a+}{$\scriptstyle{a}$}
\adjustrelabel <-.2mm, -.2mm> {b+}{$\scriptstyle{b}$}
\adjustrelabel <-.5mm, -.3mm> {c+}{$\scriptstyle{c}$}
\adjustrelabel <-.9mm, -.2mm> {d+}{$\scriptstyle{d}$}
\adjustrelabel <-.4mm, 0mm> {e+}{$\scriptstyle{e}$}
\endrelabelbox
\caption{A $w$--disc with amenable word $w =
attbt^{-1}t^{-1}cdt^{-1}et^2$ and its blown up form $t\, (s^{-1}as) \, t
\, (b) (s^{-2}cs^2) \, t^{-1} \, (s^{-1}ds)(s^{-2}es^2)$ inside} 
\endfigure

Next we replace the $w(s,t)$--disc by a smaller $w(s,t)$--disc surrounded
by ``two-leg'' discs which convert the $s$--legs to $t$--legs. A
``two-leg'' disc is labelled $h^t(h^s)\inv$ and actually has two legs
labelled $t$ but two sets of legs labelled $s$ (corresponding to the
occurrences of $s$ in $h^s$ and $h$). This is done as follows. 

Consider the $s$--legs in a disc labelled by $w(s,t)$.  They come in
sets corresponding to the $x_i(s)\in H$.  Now $x_i(s)$ lies in $J$ and
it follows that the $s$--legs can be paired off in nested cancelling
pairs.  Further by counting the occurrences of $s$ (algebraically) from
the start and considering a local minimum, we can find an adjacent
cancelling pair of the form $s\inv s$.  Deleting this pair and using
induction, we see that the pairing off can be assumed to be {\sl
coherent}; in other words cancelling pairs are always of this form
(rather than $ss\inv$).  We use
induction on the number of such pairs.  Consider an outermost pair.
What lies inside also lies in $J$ and it follows that there is two-leg
disc which converts this pair into $t$--legs.  By induction we can
find a set of two-leg discs to convert the remaining $s$--legs in
similar cancelling pairs into $t$--legs.  The effect is pictured in
figure 5.

\figuresc
\relabelbox
\epsfbox{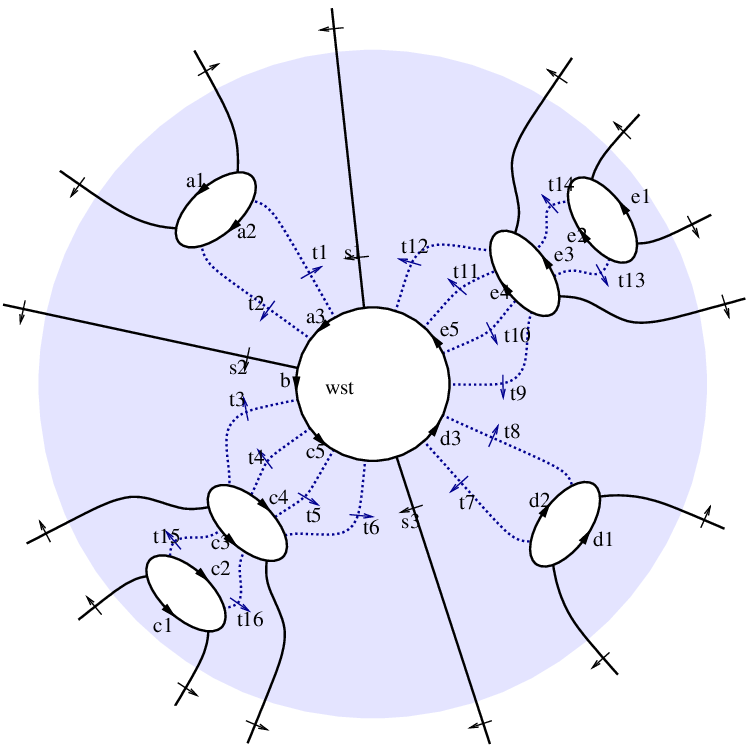}
\adjustrelabel <-.4mm, .4mm> {a1}{$\scriptstyle{a}$}
\adjustrelabel <.2mm, -.3mm> {a2}{$\scriptstyle{a}$}
\adjustrelabel <-.3mm, .3mm> {a3}{$\scriptstyle{a}$}
\adjustrelabel <-.6mm, -.3mm> {b}{$\scriptstyle{b}$}
\adjustrelabel <-.5mm, .5mm> {c1}{$\scriptstyle{c}$}
\adjustrelabel <0mm, -.5mm> {c2}{$\scriptstyle{c}$}
\adjustrelabel <-.3mm, -.5mm> {c3}{$\scriptstyle{c}$}
\adjustrelabel <-.4mm, .5mm> {c4}{$\scriptstyle{c}$}
\adjustrelabel <-.3mm, 0mm> {c5}{$\scriptstyle{c}$}
\adjustrelabel <-.3mm, -.7mm> {d1}{$\scriptstyle{d}$}
\adjustrelabel <-.5mm, .2mm> {d2}{$\scriptstyle{d}$}
\adjustrelabel <-.3mm, 0mm> {d3}{$\scriptstyle{d}$}
\adjustrelabel <-.2mm, -.2mm> {e1}{$\scriptstyle{e}$}
\adjustrelabel <-1mm, .5mm> {e2}{$\scriptstyle{e}$}
\adjustrelabel <-.3mm, 0mm> {e3}{$\scriptstyle{e}$}
\adjustrelabel <-.5mm, -.3mm> {e4}{$\scriptstyle{e}$}
\adjustrelabel <-.2mm, 0mm> {e5}{$\scriptstyle{e}$}
\adjustrelabel <0mm, -.2mm> {wst}{$w(s,t)$}
\adjustrelabel <-.5mm, .6mm> {s1}{$\scriptstyle{t}$}
\adjustrelabel <-.3mm, -.1mm> {s2}{$\scriptstyle{t}$}
\adjustrelabel <-.4mm, -.5mm> {s3}{$\scriptstyle{t}$}
\adjustrelabel <-.5mm, -.3mm> {t1}{$\scriptstyle{s}$}
\adjustrelabel <-.4mm, -.6mm> {t2}{$\scriptstyle{s}$}
\adjustrelabel <-.4mm, -.3mm> {t3}{$\scriptstyle{s}$}
\adjustrelabel <-.5mm, -.4mm> {t4}{$\scriptstyle{s}$}
\adjustrelabel <-.9mm, .4mm> {t5}{$\scriptstyle{s}$}
\adjustrelabel <-.7mm, 0mm> {t6}{$\scriptstyle{s}$}
\adjustrelabel <-.8mm, .4mm> {t7}{$\scriptstyle{s}$}
\adjustrelabel <-.3mm, .2mm> {t8}{$\scriptstyle{s}$}
\adjustrelabel <-.4mm, .2mm> {t9}{$\scriptstyle{s}$}
\adjustrelabel <-.6mm, .5mm> {t10}{$\scriptstyle{s}$}
\adjustrelabel <0mm, 0mm> {t11}{$\scriptstyle{s}$}
\adjustrelabel <-.3mm, 0mm> {t12}{$\scriptstyle{s}$}
\adjustrelabel <-3.5mm, -.4mm> {t13}{$\scriptstyle{s}$}
\adjustrelabel <-2mm, -2.3mm> {t14}{$\scriptstyle{s}$}
\adjustrelabel <2mm, 1mm> {t15}{$\scriptstyle{s}$}
\adjustrelabel <.6mm, 2.5mm> {t16}{$\scriptstyle{s}$}
\endrelabelbox
\caption{Converting $s$--legs to $t$--legs. The outer boundary reads
$w(t,t)$.}  
\endfigure

Doing this for each $x_i$ we create a partial diagram having $t$--arcs
and $s$--arcs, whose outer boundary reads $w(t,t)$. This diagram joins up
with the outer ring of figure 4. 

The last step is to delete every $s$--arc and insert ``$s$'' into the
labels at its endpoints. Here we are passing from a diagram over $G$ to a
diagram over $\Gamma$. The result is a $W\!$--diagram; see figure 6. 
\figure
\relabelbox
\epsfbox{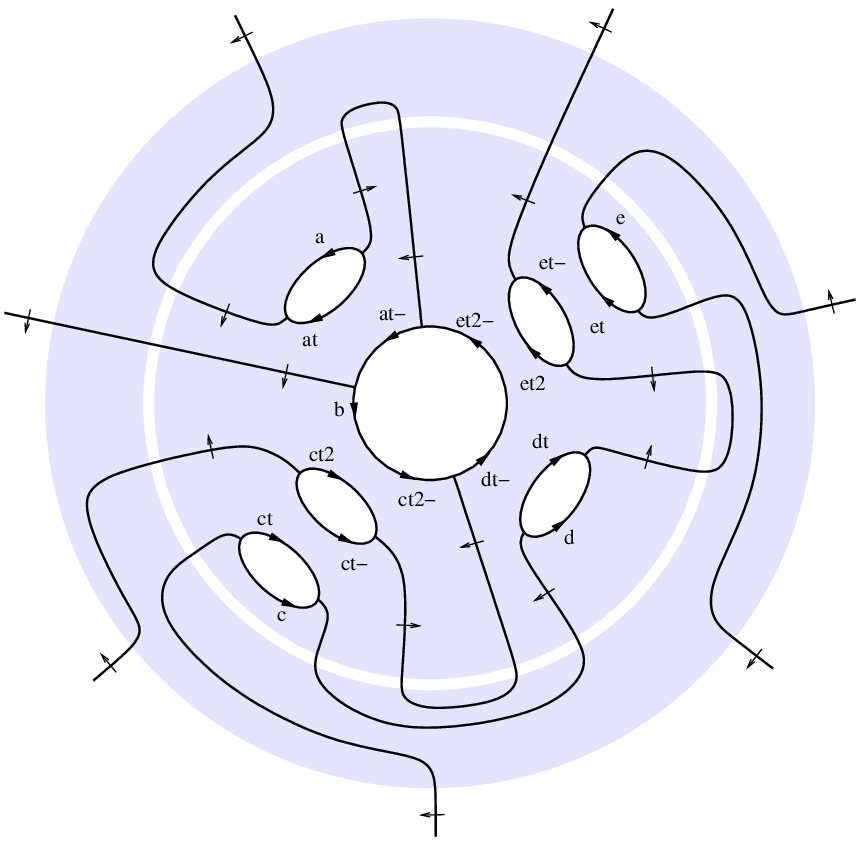}
\adjustrelabel <-.1mm, -.1mm> {a}{$\scriptstyle{a}$}
\adjustrelabel <-.2mm, 0mm> {b}{$\scriptstyle{b}$}
\adjustrelabel <-.3mm, -.3mm> {c}{$\scriptstyle{c}$}
\adjustrelabel <0mm, 0mm> {ct}{$\scriptstyle{\bar{s}cs}$}
\adjustrelabel <-1mm, 0mm> {ct-}{$\scriptstyle{\bar{s}cs}$}
\adjustrelabel <-.8mm, 0mm> {ct2}{$\scriptstyle{\bar{s}^2\!cs^{\!2}}$}
\adjustrelabel <-1.6mm, -.3mm> {ct2-}{$\scriptstyle{\bar{s}^2\!cs^{\!2}}$}
\adjustrelabel <0mm, -.3mm> {at}{$\scriptstyle{\bar{s}as}$}
\adjustrelabel <-.9mm, -.1mm> {at-}{$\scriptstyle{\bar{s}as}$}
\adjustrelabel <-.4mm, -.2mm> {d}{$\scriptstyle{d}$}
\adjustrelabel <-.2mm, .2mm> {dt}{$\scriptstyle{\bar{s}ds}$}
\adjustrelabel <-.2mm, -.2mm> {dt-}{$\scriptstyle{\bar{s}ds}$}
\adjustrelabel <.1mm, -.4mm> {e}{$\scriptstyle{e}$}
\adjustrelabel <-.7mm, .3mm> {et}{$\scriptstyle{\bar{s}es}$}
\adjustrelabel <-.3mm, -.5mm> {et-}{$\scriptstyle{\bar{s}es}$}
\adjustrelabel <-1mm, -.2mm> {et2}{$\scriptstyle{\bar{s}^2\!es^{\!2}}$}
\adjustrelabel <-.9mm, .2mm> {et2-}{$\scriptstyle{\bar{s}^2\!es^{\!2}}$}
\endrelabelbox
\caption{The $w$--disc with $t$--shape $ttt^{-1}t^{-1}t^{-1}tt$ realised
as a $W\!$--diagram with $t$--shape $ttt^{-1}$ (the root shape of
$ttt^{-1}t^{-1}t^{-1}tt$)} 
\endfigure
This $W\!$--diagram may not be irreducible because it may have strings
of two-leg discs violating condition (4). Such a string can simply be
replaced by a $t$--arc (cf [\CR; figure 4, page 134]). This has no
effect on the group relations implicit in the diagram.  This may
result in floating $t$--circles which we now delete.  The result is an
irreducible $W\!$--diagram. Finally note that every element of $Y$
(respectively $X$) is free relative to $H$ (respectively $H'$), by
[\FR, Lemma 4.3]. Thus the diagram contradicts the Root Shape
Theorem. \qed

\section{Proof of the Main Theorem}

We now return to the situation described in the Introduction. We are
given a complex $K$ with torsion-free fundamental group, and $L = K \cup
e^1 \cup e^2$, with the $2$--cell attached by an amenable word $w$. Let
$K^+ = K \cup e^1$. 

First observe that the Diagram Theorem translates directly into the
statement that the map $\pi_2(K^+,K) \to \pi_2(L,K)$ is surjective. Using
a standard transversality argument (as in [\FR, Section 3]), every
element of $\pi_2 
(L,K)$ is represented by a $w$--diagram. If the diagram contains a disc
labelled by $w$ or $\bar w$, then there must be a cancelling pair of such
discs (as in figure 2), by the Diagram Theorem. Such a pair can be
removed by a homotopy. In this way, all discs labelled by $w$ or $\bar w$
can be removed. Then the diagram represents an element of $\pi_2(K^+, K)$. 

The next lemma separates the Diagram Theorem into the adjunction problem
(condition (a)) and an additional conclusion concerning $\pi_2$. 

\proc{Lemma}\key{lemma}
Let $K \subset K^+ \subset L$ be CW complexes such that $K^+ = K \cup
\{1\hbox{--cells\/} \}$ and $L = K^+ \cup \{2\hbox{--cells\/}\}$. 
Then $\pi_2(K^+,K) \to \pi_2(L,K)$ is surjective if and only
if the following two conditions hold: 
\items
\item{\rm (a)} $\pi_1 (K) \to \pi_1 (L)$ is injective, 
\item{\rm (b)} $\pi_2 (K^+) \to \pi_2 (L)$ is surjective. 
\enditems
\endproc

\prf The Lemma is easily deduced by inspecting the braid of long exact
sequences for the triple $(L,K^+,K)$: 

$$\relabelbox
\epsfbox{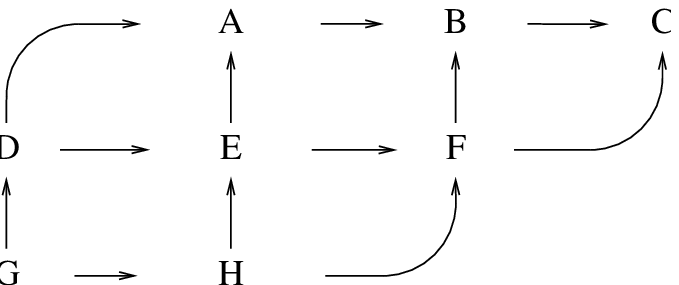}
%%% cd1.eps is exported from cd1.fig at 91.3 magnification
%%% (to undo \magstephalf ( = 1.095) in \gtmacros)
\adjustrelabel <-.68cm, 0cm> {A}{$\pi_2(L,K^+)$}
\adjustrelabel <-.5cm, 0cm> {B}{$\pi_1(K^+)$}
\adjustrelabel <-.35cm, 0cm> {C}{$\pi_1(L)$}
\adjustrelabel <-.35cm, 0cm> {D}{$\pi_2(L)$}
\adjustrelabel <-.58cm, 0cm> {E}{$\pi_2(L,K)$}
\adjustrelabel <-.4cm, 0cm> {F}{$\pi_1(K)$}
\adjustrelabel <-.52cm, 0cm> {G}{$\pi_2(K^+)$}
\adjustrelabel <-.7cm, 0cm> {H}{$\pi_2(K^+,K)$}
\endrelabelbox$$

Apart from the braid, the only information needed is that the map
$\pi_1(K) \to \pi_1(K^+)$ is injective, which is clear from the
hypotheses. 
\endprf

We also note the following basic result. 

\proc{Proposition}{\rm[\Brown, Chapter III, 5.3]}\key{brownprop}\qua 
Let $M$ be a $\Z G$--module whose underlying abelian group is a direct sum
$\bigoplus_{i} M_{i}$. Suppose that the $G$--action preserves
the direct sum decomposition and acts transitively on summands. 
Then for any $i$ we have $M \cong \Z G \otimes_{\Z G_{i}}
M_{i}$, where $G_{i}$ is the stabiliser of $M_{i}$. \qed
\endproc

\proc{Proposition}\key{mainprop}
Assume the hypotheses of Lemma \ref{lemma}, and suppose the map 
$\pi_2(K^+,K) \to \pi_2(L,K)$ is surjective. Then there is an isomorphism
of $\Z \pi_1(L)$--modules $\pi_2 (L) \cong \Z \pi_1 (L) \otimes_{\Z
\pi_1 (K)} \pi_2 (K)$. 
\endproc

This proposition (together with the Diagram Theorem) completes the proof
of the Main Theorem. 

\prf
Let $\wtilde{K}$, $\wtilde{K}^+$, and $\wtilde{L}$ be the universal
covers of $K$, $K^+$, and $L$ respectively. Let $\what{K}^+$ be the
preimage of $K^+$ in $\wtilde{L}$. Then $\wtilde{L} - \what{K}^+$
consists of $2$--cells, and there are covering maps $\wtilde{K}^+ \to
\what{K}^+ \to K^+$. Note that conditions (a) and (b) of Lemma
\ref{lemma} hold. 

By condition (a) the preimage of $K$ in $\wtilde{L}$ consists of 
disjoint copies of $\wtilde{K}$. Then $\what{K}^+$ is this disjoint
union, joined by $1$--cells. The homology group $H_2(\what{K}^+)$ is a
direct sum of copies of $H_2(\wtilde{K})$, indexed by the copies of
$\wtilde{K}$ in $\what{K}^+$. The group $\pi_1 (L)$ acts on $\what{K}^+$ by
covering translations in $\wtilde{L}$, and the induced action on
$H_2(\what{K}^+)$ preserves the direct sum decomposition (and acts
transitively on summands). 

In the action of $\pi_1 (L)$ on $\what{K}^+$, the stabilisers of the copies
of $\wtilde{K}$ are the conjugates of $\pi_1 (K)$ in $\pi_1 (L)$. These
subgroups are the stabilisers of the summands of $H_2(\what{K}^+)$ as
well. Choosing the appropriate copy of $\wtilde{K}$, Proposition
\ref{brownprop} 
implies $$H_2(\what{K}^+) \ \cong \ \Z\pi_1 (L) \otimes_{\Z \pi_1 (K)}
H_2(\wtilde{K}).$$  
So far we have used only condition (a). Now consider the map
$\pi_2 (K^+) \to \pi_2 (L)$. In the following diagram, the lower vertical
maps are Hurewicz homomorphisms, and the other maps are induced by
inclusions or coverings. 
$$ \xymatrix@R=.55cm@C=.65cm{
& \pi_2(K^+) \ar@{>>}[r] & \pi_2(L) \\
\pi_2(\wtilde{K}^+) \ar[r]^\cong \ar[d]_\cong & \pi_2(\what{K}^+) \ar[r]
\ar[d] \ar[u]^\cong & \pi_2(\wtilde{L}) \ar[d]_\cong \ar[u]^\cong \\
H_2(\wtilde{K}^+) \ar[r] & H_2(\what{K}^+) \ar[r] & H_2(\wtilde{L}) 
}$$
Condition (b) implies that the map $H_2(\what{K}^+) \to H_2(\wtilde{L})$
is surjective. It is injective because $H_3(\wtilde{L},\what{K}^+) = 0$,
and therefore there are isomorphisms $\pi_2 (L) \cong H_2 (\wtilde{L})
\cong H_2 (\what{K}^+)$. The latter is isomorphic to $\Z \pi_1 (L)
\otimes_{\Z \pi_1 (K)} \pi_2 (K)$. 
\endprf

\section{Counterexamples}

In this section we give examples showing that the Main Theorem can fail
if either condition (1) or (2) is dropped. This is in contrast with the
adjunction property which is expected to hold in much greater
generality. The adjunction property holds for the first two examples
below by [\Ho] and [\Lev]. Additional examples can be found in
[\BoPr, \Edj, \BBP, \HoMet, \Met]. 

Suppose that $L = K \cup e^1 \cup e^2$, $K^+ = K \cup e^1$, and $\pi_1(K)
\to \pi_1(L)$ is injective. Keeping the notation of the previous section, 
we have $H_2(\what{K}^+) \cong  
\Z\pi_1 (L) \otimes_{\Z \pi_1 (K)} \pi_2(K).$ We wish to construct 
examples where the injective map $H_2(\what{K}^+) \to H_2(\wtilde{L})
\cong \pi_2(L)$ is not surjective. For this it suffices to show that
$H_2(\wtilde{L}) \to H_2(\wtilde{L},\what{K}^+)$ is non-zero. 

Any $w$--diagram defines a map $S^2 \to L$ in which fibres of $t$--arcs
map to $e^1$ and discs labelled $w$ or $\bar w$ map to $e^2$. Lifting
this map to $\wtilde{L}$, one obtains a cycle representing an element of
$H_2(\wtilde{L})$. The image of this cycle in
$H_2(\wtilde{L},\what{K}^+)$ can be read off from the diagram, as a
linear combination of $2$--cells above $e^2$; see figure 7. Note that
the $2$--cells above $e^2$ form a basis for $H_2(\wtilde{L},
\what{K}^+)$, and are acted on freely and transitively by
$\pi_1(L)$.  Hence $H_2(\wtilde{L},\what{K}^+)\cong\Z\pi_1(L)$.  Given
this relative $2$--cycle, one then needs to determine whether it is
trivial, ie whether it cancels completely in $\Z
\pi_1(L)$. To decide this in general one needs to be able to solve the
word problem in $\pi_1(L)$. 

\rk{First Example} 
This is the $w$--diagram shown in figures 1 and 7. The word $w =
tatbt^{-1}c$ has amenable $t$--shape but the group $G$ has relations $b^2
= c^2 = 1$, $a^{-1}ba = c$  and so is not torsion-free. 
\figuresc
\relabelbox
\epsfbox{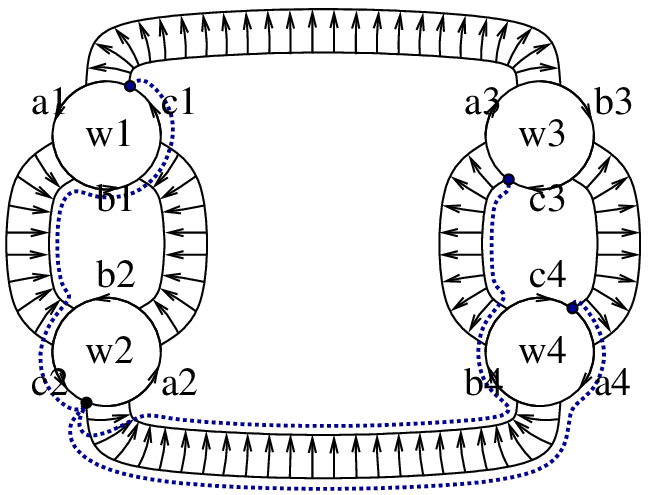}
\adjustrelabel <0mm, .2mm> {a1}{$\scriptstyle{a}$}
\adjustrelabel <-.1mm, .5mm> {a2}{$\scriptstyle{a}$}
\adjustrelabel <-.2mm, .2mm> {a3}{$\scriptstyle{a}$}
\adjustrelabel <.2mm, 0mm> {a4}{$\scriptstyle{a}$}
\adjustrelabel <.3mm, -.9mm> {b1}{$\scriptstyle{b}$}
\adjustrelabel <0mm, -.1mm> {b2}{$\scriptstyle{b}$}
\adjustrelabel <0mm, 0mm> {b3}{$\scriptstyle{b}$}
\adjustrelabel <.2mm, .2mm> {b4}{$\scriptstyle{b}$}
\adjustrelabel <.9mm, .4mm> {c1}{$\scriptstyle{c}$}
\adjustrelabel <-.1mm, -.4mm> {c2}{$\scriptstyle{c}$}
\adjustrelabel <0mm, -.1mm> {c3}{$\scriptstyle{c}$}
\adjustrelabel <0mm, 0mm> {c4}{$\scriptstyle{c}$}
\adjustrelabel <-2mm, 0mm> {w1}{$\scriptstyle{\bar{c} tb\bar{t}c D}$}
\adjustrelabel <.6mm, 0mm> {w2}{$\scriptstyle{D}$}
\adjustrelabel <-2mm, 0mm> {w3}{$\scriptstyle{-tb\bar{t}D}$}
\adjustrelabel <-1.5mm, 0mm> {w4}{$\scriptstyle{-\bar{a}\bar{t}D}$}
\endrelabelbox
\caption{The diagram as an element of $H_2(\wtilde{L},
\what{K}^+)$}
\endfigure
Choose $K$ with fundamental group $G = \la a,b,c \mid b^2 = c^2 = 1, \
a^{-1}ba = c \ra$ (the {\sl universal group} of the diagram, in other
words the group with generators the coefficients of $w$ and relations
read from the inside regions).  Note that $\pi_1(K) \to \pi_1(L)$ is
injective by [\Ho, Theorem 2], so the preceding discussion applies.
The diagram represents a relative cycle of the form $D + \bar{c}
tb\bar{t}c D - tb\bar{t}D - \bar{a}\bar{t}D$ in $H_2(\wtilde{L},
\what{K}^+)$. To see that this is non-trivial it suffices to check that
$tb\bar{t}$ and $\bar{a}\bar{t}$ are non-trivial in $\pi_1(L)$, for then
the term $D$ does not cancel with anything. For this we note that
$\pi_1(L) = \la G, t \mid tatbt^{-1}c\ra$ acts as a reflection
group in the plane. Take $\ell_1$ and $\ell_2$ to be lines meeting with
angle $\pi/3$, and let $a$, $b$, and $c$ be reflection across $\ell_1$,
and $t$ reflection across $\ell_2$. Then $\bar{a}\bar{t}$ is an order $3$
rotation and $t b \bar{t}$ a reflection, so both are non-trivial.  

Relative presentations of the form $\la G, t \mid tatbt\inv c\ra$ are
explored fully in [\Edj]. In this paper (building on results 
in [\BoPr]) a near-complete characterisation is given of groups that
admit irreducible diagrams based on $w = tatbt\inv c$. 
\endrk

\rk{Second Example}
Let $w$ have $t$--shape $t^n$ for some $n>1$.  As noted earlier this
$t$--shape is not amenable. Consider the $w$--diagram consisting of
two discs, labelled $w$ and $\bar w$, with $t$--arcs joined
cyclically. However these $t$--arcs should be joined so that the
diagram is irreducible, meaning that the two basepoints do not line
up. The group relations implied by this diagram simply say that
various coefficients are equal, in fact that $w$ itself is a proper
power. So this irreducible diagram is valid over any group $G$.  Note
that the universal group of this diagram is in fact a free group.

To see that this diagram is non-zero in $H_2(\wtilde{L},\what{K}^+)$
consider the map $G\to 1$ (the trivial group).  This is realised by a
map of $K$ to a point $K_0$ and $L$ to $L_0=e^1\cup_ne^2$ (the {\sl pseudo
projective plane} $\Bbb P_n$ with fundamental group $\Z / n\Z$). The
diagram now represents a non-trivial element of the form $D - t^kD$ in
$H_2(\wtilde{L}_0,\what{K}_0^+)$ where $k$ is the relative shift
between basepoints. In fact this diagram, regarded as a map $S^2 \to
\Bbb P_n$, represents a standard generator of $\pi_2(\Bbb P_n) \cong
\Z^{n-1}$, the others being obtained by varying $k$.

This example illustrates the difference between the notions of
irreducibility used here and in [\BoPr], where this diagram would be
considered reducible. 
\endrk

\rk{Substitution}

Before embarking on the third example, we remark that there is a very
simple way to enlarge any diagram by {\sl substitution}.  Let $w\in
G*\gen t$ be cyclically reduced and let $H$ be a group containing $G$.
Let $u$ be any reduced word in $H \ast \gen t$ which starts and
finishes with $t^{\pm1}$.  Let $w'$ be the word in $H \ast\gen t$
obtained by substituting $u^\epsilon$ for $t^\epsilon$ throughout.  We
want $w'$ to also be cyclically reduced.  In other words we require
that no cancellation between copies of $u$ is possible after
substitution.  To ensure this we need to assume that if $u$ starts and
finishes with $t^\epsilon, t^{-\epsilon}$ then {\sl all\/} the
coefficients of $w$ are non-trivial (ie middle coefficients as well as
the top and bottom ones).

If $X$ is any $w$--diagram then we can convert it to a $w'$--diagram
$X'$ by replacing the $t$--arcs by $u$--arcs.  In other words we
replace them by parallel sets of $t$--arcs corresponding to the
occurrences of $t$ in $u$.  The new regions contribute relations in $H$
of the form $xx\inv$ for the various coefficients $x$ of $u$.  If $X$
is irreducible then so is $X'$; see figure 8.

\rk{Third Example}
Let $H\supset G$ be any groups, let $w\in G*\gen t$ have $t$--shape
$t^n$ for some $n>1$ and let $u$ be any reduced word in $H *\gen t$ as
above.  Let $X$ be the $w$--diagram in the second example and $X'$
the $w'$--diagram obtained by substitution.  The $t$--shape of $w'$
is {\sl periodic} and we can clearly obtain any periodic $t$--shape in
this way.

A complete example of this construction is shown in figure 8.  This is
the simplest case, namely $w=(at)^2$ and $u=tbt\inv$.

\figuresc
\relabelbox
\epsfbox{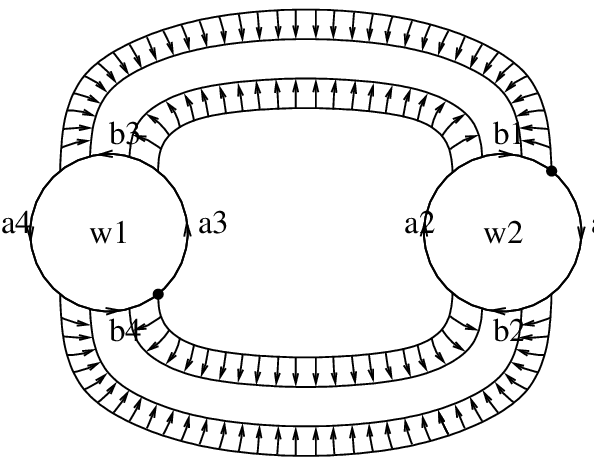}
\adjustrelabel <0mm, -.5mm> {a1}{$\scriptstyle{a}$}
\adjustrelabel <-1mm, -.5mm> {a2}{$\scriptstyle{a}$}
\adjustrelabel <0mm, -.5mm> {a3}{$\scriptstyle{a}$}
\adjustrelabel <0mm, -.5mm> {a4}{$\scriptstyle{a}$}
\adjustrelabel <0mm, 0mm> {b1}{$\scriptstyle{b}$}
\adjustrelabel <0mm, .5mm> {b2}{$\scriptstyle{b}$}
\adjustrelabel <-.5mm, 0mm> {b3}{$\scriptstyle{b}$}
\adjustrelabel <-.5mm, .5mm> {b4}{$\scriptstyle{b}$}
\adjustrelabel <.5mm, 0mm> {w1}{$w'$}
\adjustrelabel <.5mm, 0mm> {w2}{$\overline{w}'$}
\endrelabelbox
\caption{The result of substitution}
\endfigure

The universal group for $X'$ is again a free group.  We cannot prove
that $X'$ is non-zero in \strut$H_2(\wtilde{L},\what{K}^+)$ for all
choices of this construction (and we suspect that it sometimes fails
to be so) but for many choices we can prove this.  In particular, if
we take one of the coefficients $x$ of $u$ to generate a free factor
of $H=\gen x *H_0$ and map $H_0$ to the trivial group and $t$ to 1,
then this example again maps to a generator of $H_2(\Bbb P_n)$ (with
$x$ representing the 1--cell) as in the last example.

\rk{Limits to Klyachko's methods}

We have seen that Klyachko's methods imply the Diagram Theorem which
in turn proves more than is necessary for the adjunction problem for
torsion-free groups.  We have also seen that this extra information is
false for periodic $t$--shapes.  This implies that Klyachko's methods
will need modification if they are to prove the adjunction problem for
periodic $t$--shapes.  It should be noted that our examples where the
group is torsion-free are reducible if one is allowed to change the
basepoint on some $w$--discs when there is symmetry.  Call an example
in which this is not possible {\sl strongly irreducible\/}.  Thus one
necessary modification will be to incorporate an argument which makes
strong irreducibility necessary.

It is possible to construct strongly irreducible examples, with power
$t$--shape, with torsion in the universal group, which have the
property that the torsion is not implied by any single region, but is
a consequence of several.  [The simplest example has $t$--shape $t^3$
and is based in the 1--skeleton of a cube.  Labels are added to make
the diagram irreducible.  There are two ways of doing this up to
obvious symmetries.  The less symmetrical labelling is the required
diagram.]  The crash argument as presently formulated finds a single
region where the torsion-free hypothesis is contradicted.  So if the
methods are to work in general they must also be extended (perhaps by
extending the crash theorems) to provide a contradiction involving
several regions.

We have not found an example of a strongly irreducible diagram over
a torsion-free group.

Finally it should also be noted that there are $t$--shapes which are
neither periodic nor amenable.  The simplest example is
$t(tt^{-1})t(tt^{-1})^2$ and for such shapes we have no information at
all about either the adjunction problem or the Diagram Theorem.

\section{Further results} 

In [\CR] Cohen and Rourke proved that $gt$ is never in the normal
closure $\ngen w$ of $w$ in $G*\gen t$.  This result has an analogue
for diagrams which generalises the Diagram Theorem.  In [\CR; section
6] the Cohen--Rourke result is extended to prove that if $w$ has
amenable $t$--shape then no word of $t$--shape $t^n$ is in $\ngen w$.
Implicit in the proof of this is a stronger result which we shall
deduce from the Extended Diagram Theorem, see the corollary below.

Let $G$ be a group and $w\in G*\gen{t}$.  A {\sl $w$--diagram with
boundary} means a finite collection of discs in the standard 2--disc
$D^2$ labelled by $w$ and $\bar w$ together with $t$--arcs which
complete the legs compatibly with the orientations or which terminate
at the boundary.  The diagram cuts $D^2$ into a number of inside
regions (which are disjoint from the boundary $S^1$) and a finite
number of outside regions which meet the boundary.  We require that
the word in $G$ obtained by reading anticlockwise around each inside
region is the identity in $G$.  Reading round the part of the
boundary of an outside region disjoint from $S^1$, we read a word in
$G$ which we use to label the arc of intersection with $S^1$.  We can
now read a word $z\in G*\gen{t}$ from the labels on $S^1$ and the
$t$--arcs arcs which terminate on $S^1$.  We call $z$ the {\sl boundary
word}.  Irreducibility has exactly the same meaning for diagrams with
boundary as for ordinary diagrams.

Recall that the complexity of a $t$--shape is the number of Magnus
differentiations required to reduce the shape to a power shape (ie
$t^q$ for $q\in\Z$).  The {\sl complexity\/} of a word in $G*\gen{t}$
is the complexity of its $t$--shape. 

\proclaim{Extended Diagram Theorem}
Suppose that $w$ is a word in $G*\gen{t}$ with an amenable $t$--shape and
that $G$ is a torsion-free group. Then there are no irreducible
$w$--diagrams with boundary a word of complexity strictly smaller than 
the complexity of $w$. 
\endproc 

\proclaim{Corollary}{\rm (Extension of Cohen--Rourke Extension [\CR; 
page 141])}\qua
Suppose that $w$ is a word in $G*\gen{t}$ with an amenable $t$--shape
and that $G$ is a torsion-free group. Then no word of complexity
strictly smaller than the complexity of $w$ lies in the normal closure
of $w$ in $G*\gen t$.
\endproc 

Notice that this result is sharp in the sense that there are obviously
words of complexity the same as $w$ which do lie in $\ngen w$ (and
obvious corresponding diagrams with boundary).  Note also the analogy
with small cancellation theory.

The proof of the Extended Diagram Theorem is a combination of the
proof in [\CR; page 141] and the proof of the Diagram Theorem in this
paper.  We first prove an Extended Root Shape Theorem, which is a
version for $W\!$--diagrams where $W$ has the same meaning as in the
Root Shape Theorem and the boundary is a word of $t$--shape $t^q$ for
some $q\in\Z$.  To prove this we consider the cell subdivision of the
2--sphere obtained by putting a vertex in each region and joining by
edges across $t$--arcs.  There is one outside cell $C$ which has $q$
edges crossing the $t$--arcs terminating on $S^1$.  A traffic flow is
defined as before with the car on the boundary of $C$ dealt with as
described in [\CR; page 141]: all the edges of $C$ are oriented the
same way (without loss let this be ``uphill'').  Notice that any other
cell with an edge in common with $C$ has its car traverse that edge in
the ``downhill" direction, since adjacent cells induce opposite
orientations on a common edge.  Choose any point $\omega\in\d C$ not
at a vertex.  The flow constructed as before for cells other than $C$
has the property that there are times when all cars are going uphill
and hence are not on $\d C$.  This leaves time for car on $\d C$ to
rush round from just after $\omega$ to just before and hence there are
no complete crashes on $\d C$ except at $\omega$.  There must be
another crash and this leads to the identical contradiction as in the
proof of the Root Shape Theorem.

To prove the Extended Diagram Theorem we use the argument given in the
proof of the ordinary Diagram Theorem to convert any $w$-diagram to a
$W\!$--diagram: we convert the $w$--discs to partial $W\!$--diagrams
exactly as illustrated in figures 5 and 6 and we do the same for the
boundary.  Using the same notation as before, $w$ has complexity $m$.
A word of complexity less than $m$ has root $t^q$ say and is obtained
from $t^q$ by blowing up, ie, by inserting elements of $J$.  Each such
element can be converted to an element of $G*\gen s$ by using 2--leg
discs.  This converts the boundary into a word of shape $t^q$.  The
final diagram contradicts the Extended Root Shape Theorem.

\references 
\bye